\documentclass[a4paper,12pt]{article} 

\usepackage{color}

\usepackage[pagewise]{lineno}
\usepackage{amsfonts, amsmath, amsthm, amssymb}
\usepackage[T1]{fontenc}
\usepackage[cp1250]{inputenc}
\usepackage{xcolor}
\usepackage{graphicx}
\usepackage{amssymb}
\usepackage{amsmath}
\usepackage{mathptmx}
\usepackage{helvet}
\usepackage{courier}
\usepackage{txfonts}
\usepackage{tikz} 
\usetikzlibrary{arrows}
\usepackage{type1cm}

\usepackage{verbatim}

\usepackage{graphicx}
\usepackage{epsfig,amscd,amssymb,amsxtra,amsmath,amsthm}
\usepackage{type1cm}
\usepackage[T1]{fontenc}
\usepackage{graphics}
\usepackage[mathscr]{eucal}
\usepackage[all]{xy}
\usepackage{amsmath,amscd}

\newtheorem{theorem}{Theorem}[section]
\newtheorem{proposition}[theorem]{Proposition}
\newtheorem{definition}[theorem]{Definition}
\newtheorem{lemma}[theorem]{Lemma}

\newtheorem{example}[theorem]{Example}

\newtheorem{corollary}[theorem]{Corollary}

\newtheorem{observation}[theorem]{Observation}

\newtheorem*{nonumclaim}{Claim}

\newcommand{\diam} {\mathop{\rm diam}\nolimits}

\newcommand{\ent}{\textup{ent}}

\begin{document}

\def\joinrel{\mkern-3mu}
\newcommand{\varproj}{\displaystyle \lim_{\multimapinv\joinrel-\joinrel-}}

\title{Closed relations with non-zero entropy that generate no periodic points}
\author{Iztok Bani\v c,  Goran Erceg, and Judy Kennedy}
\date{}

\maketitle

\begin{abstract}
\noindent The paper is motivated by E. ~Akin's book about dynamical systems and  closed relations \cite{A}, and by J. ~Kennedy's and G.~Erceg's recent paper about the entropy of closed relations on closed intervals \cite{EK}.  

In present paper,  we introduce the entropy of a closed relation $G$  on any compact metric space $X$  and show its basic properties.  We also introduce when such a relation $G$ generates a periodic point or finitely generates a Cantor set.  Then we show that periodic points, finitely generated Cantor sets,  Mahavier products and the entropy  of closed relations are preserved by topological conjugations.  Among other things, this generalizes the well-known results about the topological conjugacy of continuous mappings.  Finally, we prove a theorem, giving sufficient conditions for a closed relation $G$ on $[0,1]$ to have a non-zero entropy. Then we present various examples of closed relations $G$ on $[0,1]$  such that 
(1) the entropy of $G$ is non-zero,
(2) no periodic point or exactly one periodic point is generated by $G$, and
(3) no Cantor set is finitely generated by $G$.
\end{abstract}
\-

\section{Introduction}\label{s0}

Our work is motivated by E. ~Akin's book "General Topology of Dynamical Systems" \cite{A}, where dynamical systems using closed relations are presented.  Our work is also motivated by J. ~Kennedy's and G.~Erceg's paper "Topological entropy on closed sets in $[0,1]^2$" \cite{EK}, where the idea of topological entropy is generalized to closed subsets of $[0,1]\times[0,1]$, and later to closed subsets of $[0,1]^n$, for $n$ a positive integer greater than $1$.  The problem of computing topological entropy was reduced to  counting the ``boxes'' (elements of the grid covers) that certain sets generated by certain closed subsets of $[0,1]^{n}$ intersect.

In that paper, many examples of closed subsets of $[0,1]\times[0,1]$ having positive topological entropy are considered. However, all the specific examples that are considered in the paper have the property that the positive entropy is ``carried''  by a finitely generated invariant Cantor set, or, sometimes, many of them. The examples also admit many periodic points. The question then arose as to whether this was always the case. That is, could one obtain a closed relation or a closed set with positive topological entropy  that did not admit periodic points nor a finitely generated invariant Cantor set having positive topological entropy? The answer for closed and connected relations on $[0,1]$ with domain and range being  $[0,1]$ is ``yes'', almost, since the example of such a function that we present in our paper does have exactly one fixed point and no other periodic points. (Note that such a relation must have a fixed point.) We also present many examples of closed relations $G$ on $[0,1]$   such that 
\begin{enumerate}
\item $G$ is not connected,
\item the entropy of $G$ is non-zero,
\item no periodic point is generated by $G$, and
\item no Cantor set is finitely generated by $G$.
\end{enumerate}

We need to develop a fair amount of machinery in order to talk about entropy in this setting. Many authors find the idea fascinating, mysterious and a little mercurial. Entropy is often quite a local property, in that small sets in a space can determine the entropy of the whole space.  The entropy of a system is a nonnegative number and $\infty$ is also allowed. But what does this number tell you? It does not give the whole story, at least that is clear. But it is measuring something.

 Suppose $X$ is a compact metric space. Recall that if $f:X \to X$ is a continuous function, the \textit{inverse limit space} generated by $f$ is 
\begin{equation*}
 \varprojlim(X,f):=\Big\{(x_{1},x_{2},x_{3},\ldots ) \in \prod_{i=1}^{\infty} X \ | \ 
\text{ for each positive integer } i,x_{i}= f(x_{i+1})\Big\},
\end{equation*}
\noindent which we can abbreviate as $\underset{\longleftarrow }{\lim }f$. The map $f$ on $X$ induces a natural homeomorphism $\sigma$ on $\underset{\longleftarrow }{\lim }f$, which is called the \textit{shift map}, and is defined by 
$$\sigma(x_{1},x_{2},x_3, x_4,\ldots)=(x_{2},x_{3},x_4,\dots)$$ for each $(x_{1},x_{2},x_3, x_4,\ldots)$ in  $\underset{\longleftarrow }{\lim }f$.

 Generalized inverse limits, or inverse limits with set-valued functions, are a generalization of (standard) inverse limits.  Here, rather than beginning with a continuous function $f$ from a compact metric space $X$ to itself, we begin with an upper semicontinuous function $f$ from $X$ to the closed subsets of $X$. The \textit{generalized inverse limit}, or the \textit{inverse limit with set-valued mappings}, associated with this set-valued function is the set 
\begin{equation*}
\varproj(X,f):=\Big\{(x_1,x_{2},x_3,\ldots )\in \prod_{i=1}^{\infty} X \ | \  
\text{ for each positive integer }i,x_{i}\in f(x_{i+1})\Big\},
\end{equation*}
\noindent which is a closed subspace of   $\Pi _{i =1}^{\infty}X$  endowed with the product topology. (As is the case with standard inverse limits, these can be defined in much more general settings, but we do not need those here.) Here again, the shift map $\sigma$ defined above takes $\varproj(X,f)$ onto itself, but it is no longer a homeomorphism: $\sigma: \varproj(X,f) \to \varproj(X,f)$ is a continuous surjection.  And again, we often abbreviate $\varproj(X,f)$ as $\varproj f$. The topic of generalized inverse limits  is currently an
intensely studied area of continuum theory, with papers from many authors at this point. The topic of entropy is another topic currently being studied by many authors. It continues to be defined in new settings.

Generalized inverse limits, or inverse limits with set-valued functions, a subject studied only since 2003 with its introduction by Bill Mahavier, and much subsequent development by Tom Ingram, provide  an entirely new way to study multi-valued functions, a way that does not lose information under iteration. (For the interested reader, we recommend the following references: 
\cite{B}, \cite{BCMM1}, \cite{BCMM2}, \cite{BCMM3}, \cite{BK}, \cite{CR} \cite{GK1}, \cite {GK2}, \cite{Il}, \cite{IM2}, \cite{I1}, \cite{I2}, \cite{I3}, \cite{I4}, \cite{I7},\cite{L}, \cite{M}, \cite{N1},  \cite{N2},  \cite{N3},  \cite{N4}, and \cite{V}. This list is far from exhaustive. Tom Ingram and Bill Mahavier included a chapter on these spaces in their
 book \cite{IM1}, and since then Tom Ingram has written another book on the topic, \cite{I6}, where more references may be found.)
But it is increasingly apparent that generalized inverse limits also offer a rich source of new examples of dynamical systems and continua. In fact, they offer a sort of lab in which one can make mathematical experiments - and then have a real chance, with some effort, of understanding (via some sort of coding) deeply the resulting topology and dynamics of the example, and of how the topology and dynamics are interacting. Note that the dynamics of upper semicontinuous set-valued functions are often determined by closed subsets of their graphs and this is the reason why in this paper, we mainly concentrate on dynamics of closed relations on compact metric spaces rather than on dynamics of set-valued fuctions.

We generalize the notion of topological entropy of closed relations on $[0,1]$ (i.e.,  closed subsets of $[0,1]\times [0,1]$) that is introduced in \cite{EK} to topological entropy of closed relations on arbitrary compact metric spaces $X$ (i.e., closed subsets of $X\times X$); see Section \ref{s1}. In Section \ref{s2}, we introduce the notion of periodic points that are generated by closed relations on any compact metric space and define when a space is finitely generated by such a closed relation. We also  prove that if a Cantor set is finitely generated by a closed relation on a compact metric space, then there are at least two periodic points that are generated by the same relation. In Section \ref{s3}, we generalize the notion of topological conjugacy for continuous functions to the topological conjugacy for closed relations $G$ and $H$ on compact metric spaces $X$ and $Y$, respectively, and show that periodic points, finitely generated Cantor sets, Mahavier products and the entropy of closed relations are being preserved by  topological conjugation. In Section \ref{s4}, we present and prove a theorem that gives a large family of closed relations on $[0,1]$ all having a non-zero entropy. We use this theorem in Section \ref{s5}, where we present various examples of closed relations $G$ on $[0,1]$  such that 
\begin{enumerate}
\item the entropy of $G$ is non-zero,
\item no periodic point or exactly one periodic point is generated by $G$, and
\item no Cantor set is finitely generated by $G$.
\end{enumerate}

\section{Definitions and notation} \label{s11}

In this section, basic definitions and well-known results that are needed later in the paper are presented.
\begin{definition}
Let $X$ and $Y$ be metric spaces, and let $f:X\rightarrow Y$ be a function.  We use  
$$
\Gamma(f)=\{(x,y)\in X\times Y \ | \ y=f(x)\}
$$
to denote \emph{ \color{blue} the graph of the function $f$}.
\end{definition}
\begin{definition}
If $X$ is a compact metric space, then \emph{ \color{blue}  $2^X$ }denotes the set of all  nonempty closed subsets of $X$.  
\end{definition}
\begin{definition}
Let $X$ and $Y$ be compact metric spaces. 
A function $F: X\rightarrow 2^Y$ is called \emph{ \color{blue}  a set-valued function} from $X$ to $Y$. We denote set-valued functions $F: X\rightarrow 2^Y$ by \emph{ \color{blue}  $F: X\multimap Y$}.
\end{definition}
\begin{definition}
A set-valued function  $F : X\multimap Y$ is  \emph{ \color{blue}   upper semicontinuous  at a point $x_0\in X$},  if for each
open set  $U\subseteq Y$ such that $F(x_0)\subseteq U$, there is an open set $V$ in $X$ such that
\begin{enumerate}
\item $x_0\in V$ and
\item for each $x\in V$, $F(x)\subseteq U$. 
\end{enumerate}  
The set-valued function  $F : X\multimap Y$ is  \emph{ \color{blue}   upper semicontinuous},  if it is upper semicontinuous at any point $x\in X$.
\end{definition}
\begin{definition}
 \emph{ \color{blue} The  graph $\Gamma(F)$ of a set-valued function} $F:X\multimap Y$ is the set of
all points  $(x,y)\in X\times Y$ such that $y \in F(x)$.  The set-valued function $F$ is surjective, if $\bigcup_{x\in X}F(x)=Y$. 
\end{definition}
There is a simple characterization of upper semicontinuous set-valued functions
(\cite[Proposition 11, p.\ 128]{A} and \cite[Theorem 1.2, p.\ 3]{I6}):

\begin{theorem}
\label{th:grafi}  Let $X$ and $Y$ be compact metric spaces and $F:X\multimap Y$ a set-valued function. Then $F$ is upper semicontinuous if and only if its
graph $\Gamma(F)$ is closed in  $X\times Y$. 
\end{theorem}

\begin{definition}
Let $X$ be a compact metric space and let $G\subseteq X\times X$ be a relation on $X$. If $G$ is closed in $X\times X$, then we say that $G$ is a closed relation on $X$.  
\end{definition}

\begin{definition}
Let $X$  be a set and let $G$ be a relation on $X$.  Then we define  
$$
G^{-1}=\{(y,x)\in X\times X \ | \ (x,y)\in G\}
$$
to be the inverse relation of the relation $G$ on $X$.
\end{definition}
\begin{definition}
Let $X$ be a compact metric space and let $G$ be a closed relation on $X$. Then we call
$$
\star_{i=1}^{m}G^{-1}=\Big\{(x_1,x_2,x_3,\ldots ,x_{m+1})\in \prod_{i=1}^{m+1}X \ | \ \textup{ for each } i\in \{1,2,3,\ldots ,m\}, (x_{i+1},x_i)\in G\Big\}
$$
for each positive integer $m$, the $m$-th Mahavier product of $G$, and
$$
\star_{i=1}^{\infty}G^{-1}=\Big\{(x_1,x_2,x_3,\ldots )\in \prod_{i=1}^{\infty}X \ | \ \textup{ for each positive integer } i, (x_{i+1},x_i)\in G\Big\}
$$
the infinite  Mahavier product of $G$.
\end{definition}
\begin{observation}
Let $X$ be a compact metric space, let $f:X\rightarrow X$ be a continuous function. 
Then 
$$
\star_{n=1}^{\infty}\Gamma(f)^{-1}=\varprojlim(X,f).
$$
Also, if $F:X\multimap X$ is an upper semi-continuous function, then  
$$
\star_{n=1}^{\infty}\Gamma(F)^{-1}=\varproj(X,f).
$$
\end{observation}
\begin{definition}
Let $X$ be a compact metric space and let $G$ be a closed relation on $X$.  The function  
$$
\sigma : \star_{n=1}^{\infty}G^{-1} \rightarrow \star_{n=1}^{\infty}G^{-1},
$$
 defined by 
$$
\sigma (x_1,x_2,x_3,x_4,\ldots)=(x_2,x_3,x_4,\ldots)
$$
for each $(x_1,x_2,x_3,x_4,\ldots)\in \star_{n=1}^{\infty}G^{-1}$, 
is called \emph{ \color{blue}  the shift map on $\star_{n=1}^{\infty}G^{-1}$}.  
\end{definition}

\section{Topological entropy of closed relations on compact metric spaces}\label{s1}
In present section,  we generalize the notion of topological entropy of closed relations on $[0,1]$ (i.e., topological entropy of closed subsets of $[0,1]\times [0,1]$) that is introduced in \cite{EK} to topological entropy of closed relations on compact metric spaces.
We start with the following simple definitions that we use later.
\begin{definition}
Let $X$ be a compact metric space and let $\mathcal S$ be a family of subsets of $X$. We use  \emph{ \color{blue}  $|\mathcal S|$} to denote the cardinality of $\mathcal S$.
\end{definition}
\begin{definition}
Let $X$ be a compact metric space and let $\mathcal S$ be a family of subsets of $X$. For each positive integer $n$, we use  $\mathcal S^n$ to denote the family
$$
\mathcal S^n=\{S_1\times S_2\times S_3\times \ldots \times S_n \ | \ S_1,S_2,S_3,\ldots,S_n\in \mathcal S\}.
$$ 
We call the elements $S_1\times S_2\times S_3\times \ldots \times S_n$ of $\mathcal S$ \emph{ \color{blue}  the $n$-boxes (generated by the family $\mathcal S$)}.
\end{definition}
\begin{definition}
Let $X$ be a compact metric space and let $\mathcal U$ be an open cover for $X$. We use $N(\mathcal U)$ to denote
$$
N(\mathcal U)=\min\{|\mathcal V| \ | \ \mathcal V \text{ is a finite subcover of } \mathcal U\}.
$$
\end{definition}
\begin{definition}
Let $X$ be a compact metric space, let $K$ be a closed subset of the product $\prod_{i=1}^{n}X$, and let $\mathcal U$ be a family of open subsets of $\prod_{i=1}^{n}X$ such that $K\subseteq \bigcup \mathcal U$. Then we use $N(K,\mathcal U)$ to denote
$$
N(K,\mathcal U)=\min\Big\{|\mathcal V| \ | \ \mathcal V \text{ is a subfamily of } \mathcal U \text{ such that }  K\subseteq \bigcup \mathcal V\Big\}.
$$ 
\end{definition}
\begin{observation}
Let $X$ be a compact metric space, let $K$ be a non-empty closed subset of the product $\prod_{i=1}^{n}X$, and let $\mathcal U$ be a family of open subsets of $\prod_{i=1}^{n}X$ such that $K\subseteq \bigcup \mathcal U$.  Note that since $K$ is a non-empty compactum, $N(K,\mathcal U)$ is a positive integer.
\end{observation}

\begin{observation}\label{obs1}
Let $X$ be a compact metric space, let $G$ be a closed relation on $X$, and let $\alpha $ be an open cover for $X$.  Note that
$$
N(\star_{i=1}^{m}G^{-1},\alpha^{m+1})\leq N(\alpha )^{m+1}
$$
for each positive integer $m$.
\end{observation}

\begin{lemma}\label{the lemma we need}
Let $X$ be a compact metric space, let $G$ be a non-empty closed relation on $X$, and let $\alpha $ be an open cover for $X$. Then
$$
N(\star_{i=1}^{m+n}G^{-1},\alpha^{m+n+1})\leq N(\star_{i=1}^{m}G^{-1},\alpha^{m+1})\cdot N(\star_{i=1}^{n}G^{-1},\alpha^{n+1})
$$
for all positive integers $m$ and $n$.
\end{lemma}
\begin{proof}
Let  $m$ and $n$ be positive integers and let
\begin{enumerate}
\item $N(\star_{i=1}^{m+n}G^{-1},\alpha^{m+n+1})=k_{m+n}$ and let $\beta=\{\textup{B}_1,\textup{B}_2,\textup{B}_3,\ldots ,\textup{B}_{k_{m+n}}\}$  be a subfamily of $\alpha^{m+n+1}$ such that $\star_{i=1}^{m+n}G^{-1}\subseteq \textup{B}_1\cup \textup{B}_2\cup \textup{B}_3\cup \ldots \cup \textup{B}_{k_{m+n}}$,
\item $N(\star_{i=1}^{m}G^{-1},\alpha^{m+1})=k_m$ and let $\gamma=\{\Gamma_1,\Gamma_2,\Gamma_3,\ldots ,\Gamma_{k_{m}}\}$  be a subfamily of $\alpha^{m+1}$ such that $\star_{i=1}^{m}G^{-1}\subseteq \Gamma_1\cup \Gamma_2\cup \Gamma_3\cup \ldots \cup \Gamma_{k_{m}}$, and 
\item $N(\star_{i=1}^{n}G^{-1},\alpha^{n+1})=k_n$ and let $\delta=\{\Delta_1,\Delta_2,\Delta_3,\ldots ,\Delta_{k_{n}}\}$  be a subfamily of $\alpha^{n+1}$ such that $\star_{i=1}^{n}G^{-1}\subseteq \Delta_1\cup \Delta_2\cup \Delta_3\cup \ldots \cup \Delta_{k_{n}}$. 
\end{enumerate} 
For each $i\in \{1,2,3,\ldots ,k_m\}$,  let $A_1^{i,\Gamma}$, $A_2^{i,\Gamma}$, $A_3^{i,\Gamma}$, $\ldots$, $A_m^{i,\Gamma}$, $A_{m+1}^{i,\Gamma}\in \alpha$ such that 
$$
\Gamma_i=A_1^{i,\Gamma}\times A_2^{i,\Gamma}\times A_3^{i,\Gamma}\times \ldots \times A_m^{i,\Gamma}\times A_{m+1}^{i,\Gamma}
$$
and for each $j\in \{1,2,3,\ldots ,k_n\}$,  let $A_1^{j,\Delta}$, $A_2^{j,\Delta}$, $A_3^{j,\Delta}$, $\ldots$  ,$A_n^{j,\Delta}$, $A_{n+1}^{j,\Delta}\in \alpha$ such that 
$$
\Delta_j=A_1^{j,\Delta}\times A_2^{j,\Delta}\times A_3^{j,\Delta}\times \ldots \times A_n^{j,\Delta}\times  A_{n+1}^{j,\Delta}.
$$
Let 
$$
\mathbf x=(x_1,x_2,x_3,\ldots,x_m,x_{m+1},\ldots  ,x_{m+n},x_{m+n+1})\in \star_{i=1}^{m+n}G^{-1}
$$ 
be any point.  Then 
$$
(x_1,x_2,x_3,\ldots,x_m,x_{m+1})\in  \star_{i=1}^{m}G^{-1}
$$
and
$$
(x_{m+1},x_{m+2},x_{m+3},\ldots,x_{m+n},x_{m+n+1})\in  \star_{i=1}^{n}G^{-1}.
$$
Let $i\in \{1,2,3,\ldots ,k_m\}$ and $j\in \{1,2,3,\ldots ,k_n\}$ be such that 
$(x_1,x_2,x_3,\ldots,x_m,x_{m+1})\in\Gamma_i$ and $(x_{m+1},x_{m+2},x_{m+3},\ldots,x_{m+n},x_{m+n+1})\in \Delta_j$. Then 
$$
\mathbf x\in A_1^{i,\Gamma}\times A_2^{i,\Gamma}\times A_3^{i,\Gamma}\times \ldots \times A_m^{i,\Gamma}\times (A_{m+1}^{i,\Gamma}\cap A_1^{j,\Delta})\times A_2^{j,\Delta}\times A_3^{j,\Delta}\times \ldots \times A_n^{j,\Delta}\times  A_{n+1}^{j,\Delta}.
$$
Since 
$$
A_1^{i,\Gamma}\times A_2^{i,\Gamma}\times A_3^{i,\Gamma}\times \ldots \times A_m^{i,\Gamma}\times (A_{m+1}^{i,\Gamma}\cap A_1^{j,\Delta})\times A_2^{j,\Delta}\times A_3^{j,\Delta}\times \ldots \times A_n^{j,\Delta}\times  A_{n+1}^{j,\Delta}
$$
$$
\subseteq A_1^{i,\Gamma}\times A_2^{i,\Gamma}\times A_3^{i,\Gamma}\times \ldots \times A_m^{i,\Gamma}\times A_1^{j,\Delta}\times A_2^{j,\Delta}\times A_3^{j,\Delta}\times \ldots \times A_n^{j,\Delta}\times  A_{n+1}^{j,\Delta},
$$
it follows that $\star_{i=1}^{m+n}G^{-1}$ can be covered by $ k_m\cdot k_n$ sets from $\alpha^{m+n+1}$. Therefore, $k_{m+n}\leq k_m\cdot k_n$. 
\end{proof}
We use the following lemma in the proof of Theorem \ref{limit exists}.
\begin{lemma}\label{the sequence}
Let $(a_m)$ be a sequence in $\mathbb R$ such that
\begin{enumerate}
\item for each positive integer $m$, $a_m\geq 0$, and
\item for all positive integers $m$ and $n$, 
$$
a_{m+n}\leq a_m+a_n.
$$
\end{enumerate}
Then the limit 
$$
\lim_{m\to \infty}\frac{a_m}{m}
$$
exists and it equals to
$$
\lim_{m\to \infty}\frac{a_m}{m}=\inf\Big\{\frac{a_m}{m} \ | \ m \textup{ is a positive integer}\Big\}.
$$
\end{lemma}
\begin{proof}
The proof can be found in \cite[Theorem 4.9, page 87]{Walters book}.
\end{proof}
\begin{theorem}\label{limit exists}
Let $X$ be a compact metric space, let $G$ be a non-empty closed relation on $X$, and let $\alpha $ be an open cover for $X$.  Then the limit
$$
\lim_{m\to \infty}\frac{\log N(\star_{i=1}^{m}G^{-1},\alpha^{m+1})}{m}
$$
exists.
\end{theorem}
\begin{proof}
For each positive integer $m$,  let 
$$
a_m=\log N(\star_{i=1}^{m}G^{-1},\alpha^{m+1}).
$$
Since $G\neq \emptyset$, it follows that for each positive integer $m$, $ a_m\geq 0$. Let $m$ and $n$ be any positive integers.  By Lemma \ref{the sequence} it suffices to show that  
$$
a_{m+n}\leq a_m+a_n. 
$$
It follows (using Lemma \ref{the lemma we need}) that
$$
a_{m+n}=\log N(\star_{i=1}^{m+n}G^{-1},\alpha^{m+n+1})\leq \log(N(\star_{i=1}^{m}G^{-1},\alpha^{m+1})\cdot N(\star_{i=1}^{n}G^{-1},\alpha^{n+1}))
$$
$$
=\log N(\star_{i=1}^{m}G^{-1},\alpha^{m+1})+\log N(\star_{i=1}^{n}G^{-1},\alpha^{n+1})=a_m+a_n.
$$
\end{proof}
\begin{definition}
Let $X$ be a compact metric space, let $G$ be a non-empty closed relation on $X$, and let $\alpha $ be an open cover for $X$.  We define \emph{ \color{blue}  the entropy of $G$ with respect to the open cover $\alpha $} by
$$
\ent(G,\alpha )=\lim_{m\to \infty}\frac{\log N(\star_{i=1}^{m}G^{-1},\alpha^{m+1})}{m}.
$$
\end{definition}

\begin{definition}
Let $X$ be a metric space and let $\mathcal S$ and $\mathcal T$ be families of subsets of $X$. We say that the family $\mathcal S$ \emph{ \color{blue}  refines} the family $\mathcal T$, if for each $S\in \mathcal S$ there is $T\in \mathcal T$ such that $S\subseteq T$.  The notation 
$$
\mathcal T\leq \mathcal S
$$
means that the family $\mathcal S$ refines the family $\mathcal T$.
\end{definition}
\begin{proposition}\label{prop1}
Let $X$ be a compact metric space and let $G$ be a non-empty closed relation on $X$. For all open covers $\alpha$ and $\beta$, the following holds:
$$
\alpha\leq \beta \Longrightarrow \ent(G,\alpha )\leq \ent(G,\beta ).
$$
\end{proposition}
\begin{proof}
Let $\alpha$ and $\beta$ be any open covers for $X$ such that $\alpha\leq \beta $. Then $\alpha^n\leq \beta^n$ for each positive integer $n$.  Let $m$ be a positive integer, let $k=N(\star_{i=1}^{m}G^{-1},\beta^{m+1})$ and let
$$
\{B_1,B_2,B_3,\ldots ,B_k\}\subseteq \beta^{m+1}
$$
such that $\star_{i=1}^{m}G^{-1}\subseteq B_1\cup B_2\cup B_3\cup \ldots \cup B_k$. For each $i\in\{1,2,3,\ldots ,k\}$, let $A_i\in \alpha$ such that $B_i\subseteq A_i$. Therefore, $\star_{i=1}^{m}G^{-1}\subseteq A_1\cup A_2\cup A_3\cup \ldots \cup A_k$ and 
$$
N(\star_{i=1}^{m}G^{-1},\alpha^{m+1})\leq N(\star_{i=1}^{m}G^{-1},\beta^{m+1})
$$
 follows.  Therefore,  $\ent(G,\alpha )\leq \ent(G,\beta )$.
\end{proof}
\begin{proposition}\label{prop2}
Let $X$ be a compact metric space and let $\alpha$ be an open cover for $X$.  For all non-empty closed relations $H$ and $G$  on $X$ the following holds:
$$
H\subseteq G \Longrightarrow \ent(H,\alpha )\leq \ent(G,\alpha ).
$$
\end{proposition}
\begin{proof}
The proposition follows from the fact that 
$$
N(\star_{i=1}^{m}H^{-1},\alpha^{m+1})\leq N(\star_{i=1}^{m}G^{-1},\alpha^{m+1})
$$
for each positive integer $m$.
\end{proof}
\begin{definition}
Let $X$ be a compact metric space, let $G$ be a non-empty closed relation on $X$, and let 
$$
E=\{\ent(G,\alpha ) \ | \ \alpha \textup{ is an open cover for } X\}.
$$
  We define \emph{ \color{blue}  the entropy of $G$} by
$$
\ent(G)=\begin{cases}
				0\text{;} & G=\emptyset \\
				\sup (E)\text{;} & G\neq \emptyset \textup{ and } E \textup{ is bounded in } \mathbb R \\
				\infty\text{;} & G\neq \emptyset \textup{ and } E \textup{ is not bounded in } \mathbb R.
			\end{cases}
$$
\end{definition}
Note that this  generalizes the entropy that is defined in \cite{EK} from closed subsets of $[0,1]\times [0,1]$ to closed subsets of $X\times X$ for any compact metric space $X$. 


\begin{theorem}\label{thm2thm}
Let $X$ be a compact metric space.  For all non-empty closed relations $H$ and $G$ on  $X$, the following holds:
$$
H\subseteq G \Longrightarrow \ent(H)\leq \ent(G).
$$
\end{theorem}
\begin{proof}
The proposition follows directly from Proposition \ref{prop2}.
\end{proof}
\begin{theorem}\label{thm3thm}
Let $X$ be a compact metric space  and let $G$ be a closed relation on $X$.  Then 
$$
\ent(G^{-1}) = \ent(G).
$$
\end{theorem}
\begin{proof}
Let $m$ be a positive integer. Note that 
$$
(x_1,x_2,x_3,\ldots ,x_m,x_{m+1})\in \star_{i=1}^{m}G^{-1} \Longleftrightarrow (x_{m+1},x_m,x_{m-1},\ldots ,x_2,x_{1})\in \star_{i=1}^{m}G.
$$
Let $\alpha$ be an open cover for $X$, let $k=N(\star_{i=1}^{m}G^{-1},\alpha^{m+1})$, and let 
$$
\{U_1,U_2,U_3,\ldots, U_k\}\subseteq \alpha^{m+1}
$$  
such that $\star_{i=1}^{m}G^{-1}\subseteq U_1\cup U_2 \cup U_3\cup \ldots \cup  U_k$.
For each $i\in\{1,2,3,\ldots ,k\}$, let $A_{1}^{i}$, $ A_{2}^{i}$, $ A_{3}^{i}$, $ \ldots$, $A_{m+1}^{i}\in \alpha $ such that
$$
U_i=A_{1}^{i}\times A_{2}^{i}\times A_{3}^{i}\times \ldots \times A_{m+1}^{i}
$$
and let
$$
V_i=A_{m+1}^{i}\times A_{m}^{i}\times A_{m-1}^{i}\times \ldots \times A_{1}^{i}.
$$
Then $\star_{i=1}^{m}G\subseteq V_1\cup V_2 \cup V_3\cup \ldots \cup  V_k$. Therefore,
$$
N(\star_{i=1}^{m}G,\alpha^{m+1})\leq N(\star_{i=1}^{m}G^{-1},\alpha^{m+1}).
$$
A similar argument gives 
$$
N(\star_{i=1}^{m}G,\alpha^{m+1})\geq N(\star_{i=1}^{m}G^{-1},\alpha^{m+1}).
$$
Therefore, for each open cover $\alpha$ for $X$,
$$
N(\star_{i=1}^{m}G,\alpha^{m+1})= N(\star_{i=1}^{m}G^{-1},\alpha^{m+1}).
$$
and the result follows.
\end{proof}
In Theorem \ref{povezava}, we show that the entropy of closed relations on $X$ is a generalization of the well-known topological entropy of  continuous functions $f:X\rightarrow X$.  Before that we give the following definitions and prove some auxiliary results.

\begin{definition}
Let $X$ be a set,  let $f:X\rightarrow X$ be a function and let $\mathcal S$ be a family of subsets of $X$. Then we define
$$
f^{-1}(\mathcal S)=\{f^{-1}(S) \ | \ S\in \mathcal S\}.
$$
\end{definition}

\begin{definition}
Let $X$ be a set and let $\mathcal A_1$, $\mathcal A_2$, $\mathcal A_3$, $\ldots$, $\mathcal A_m$ be families of subsets of $X$.  Then we define
$$
\vee_{i=1}^{m}\mathcal A_i=\{A_1\cap A_2\cap A_3\cap \ldots \cap A_m \ | \ \textup{ for each } i\in \{1,2,3,\ldots ,m\}, A_i\in \mathcal A_i\}.
$$
\end{definition}
\begin{definition}
Let $X$ be a compact metric space and let $f:X\rightarrow X$ be a continuous function. For any open cover $\alpha$ for $X$, we define
$$
h(f,\alpha)=\lim_{m\to\infty}\frac{\log(N(\vee_{i=0}^{m}f^{-i}(\alpha)))}{m}.
$$
\end{definition}

\begin{definition}
Let $X$ be a compact metric space and let $f:X\rightarrow X$ be a continuous function.  Also, let $E=\{h(f,\alpha) \ | \ \alpha \textup{ is an open cover for } X\}$. Then we define
$$
h(f)=\begin{cases}
				\sup (E)\text{;} & E \textup{ is bounded in } \mathbb R \\
				\infty\text{;} & E \textup{ is not bounded in } \mathbb R
			\end{cases}
$$
to be \emph{ \color{blue}  the entropy of the function $f$}.
\end{definition}
For a more thorough discussion about the topological entropy $h(f)$, see \cite{Walters book}.  
The following theorem is a well-known result about $h(f)$ that we use later.

\begin{theorem}\label{entropyshift}
Let $X$ be a compact metric space, let $f:X\rightarrow X$ be a continuous function, and let $\sigma : \star_{n=1}^{\infty}\Gamma(f)^{-1} \rightarrow \star_{n=1}^{\infty}\Gamma(f)^{-1}$ be the shift function.  Then 
$$
h(f)=h(\sigma).
$$
\end{theorem}
\begin{proof}
The proof can be found in \cite[Proposition 5.2]{Bo2}.
\end{proof}

The following theorem shows that the same property is established also for the entropy of closed subsets of $X\times X$.
\begin{theorem}\label{the same}
Let $X$ be a compact metric space, let $G$ be a closed relation on $X$ such that $p_1(G)\subseteq p_2(G)$ and let $\sigma$ be the shift map on $\star_{i=1}^{\infty}G^{-1}$. Then 
$$
\ent(G)=h(\sigma).
$$
\end{theorem}
\begin{proof}
The proof of the lemma is essentially the same as the proof of \cite[Theorem 4.3, page 121]{EK}. Therefore, we omit it and leave the details to the reader.
\end{proof}
The following theorem shows that the topological entropy of a continuous function $f:X\rightarrow X$ is just a special case of the entropy of closed subsets of $X\times X$. Therefore, the concept $\ent(G)$ is a generalization of the concept of $h(f)$.
\begin{theorem}\label{povezava}
Let $X$ be a compact metric space, and let $f:X\rightarrow X$ be a continuous function. Then
$$
\ent(\Gamma(f))=h(f).
$$
\end{theorem}
\begin{proof}
Let $\sigma$ be the shift map on $\star_{i=1}^{\infty}{\Gamma(f)}^{-1}$.  Note that $\star_{i=1}^{\infty}{\Gamma(f)}^{-1}=\underset{\longleftarrow }{\lim }(X,f)$.  By Theorem \ref{the same}, $\ent(\Gamma(f))=h(\sigma)$ and by Theorem \ref{entropyshift}, $h(f)=h(\sigma)$. Therefore, $\ent(\Gamma(f))=h(f)$. 
\end{proof}

\section{Periodic points that are generated by closed relations on compact metric spaces}\label{s2}
Here we introduce the notion of periodic points that are generated by closed relations on compact metric spaces. We also define when a space is finitely generated by such a relation and prove that if a Cantor set is finitely generated by a closed relation on a compact metric space, then there are at least two periodic points that are generated by the same relation. 

\begin{definition}
Let $X$ be a compact metric space and let $G$ be a closed relation on $X$.   We say that \emph{ \color{blue}  a periodic point is generated by $G$}, if there is a point $\mathbf x\in  \star_{n=1}^{\infty}G^{-1} $ such that for some positive integer $n$, 
$$
\sigma^n (\mathbf x)=\mathbf x.
$$
 If 
 $$
 \sigma^n (\mathbf x)\neq \mathbf x
 $$
  for each $\mathbf x\in  \star_{n=1}^{\infty}G^{-1} $ and for each positive integer $n$, then we say that \emph{ \color{blue}  no periodic point is generated by $G$}.
\end{definition}

\begin{definition}
Let $X$ be a compact metric space and let $G$ be a closed relation on $X$ and let $A$ be a subspace of $\star_{n=1}^{\infty}G^{-1} $. We say that \emph{ \color{blue}  $A$ is finitely generated by $G$}, if there is a finite collection $F=\{x_1,x_2,x_3,\ldots,x_n\}\subseteq G$ such that $\star_{n=1}^{\infty}F^{-1} =A$.  
\end{definition}
\begin{definition}
Let $X$ be a compact metric space, let $G$ be a closed relation on $X$ and let $Y$ be a compact metric space. We say that  \emph{ \color{blue}  $Y$  is finitely generated by $G$}, if there is a subspace $A$ of  $\star_{n=1}^{\infty}G^{-1} $ such that
\begin{enumerate}
\item $A$ is homeomorphic to $Y$ and
\item  $A$ is finitely generated by $G$.
\end{enumerate}
\end{definition}
\begin{theorem}\label{gaga}
Let $X$ be a compact metric space and let  $G$ be a closed relation on $X$.  If a Cantor set is finitely generated by $G$, then at least two periodic points are generated by $G$. 
\end{theorem}
\begin{proof}
Suppose that $C$ is a Cantor set in $\star_{n=1}^{\infty}G^{-1}$ and that $F=\{c_1,c_2,c_3,\ldots,c_n\}\subseteq G$ such that $\star_{n=1}^{\infty}F^{-1}=C$.  Let $\mathbf x=(x_1,x_2,x_3,\ldots )\in \star_{n=1}^{\infty}F^{-1}$ be any point.  Since $(x_{k+1},x_k)\in F$ for each positive integer $k$ and since $F$ is finite, there are positive integers $k_1$ and $k_2$ such that $k_1< k_2$ and 
$$
x_{k_1}=x_{k_2}.
$$
Choose and fix such positive integers $k_1$ and $k_2$.  Let
$$
\mathbf y=(x_{k_1}, x_{k_1+1}, \ldots ,x_{k_2-1}, x_{k_1},  x_{k_1+1}, \ldots ,x_{k_2-1}, x_{k_1},  x_{k_1+1}, \ldots ,x_{k_2-1}, x_{k_1}, \ldots).
$$
Then $\mathbf y\in \star_{n=1}^{\infty}G^{-1}$  and $\sigma^{k_2-k_1}(\mathbf y)=\mathbf y$ which means that $\mathbf y$ is a periodic point that is generated by $G$.  Obviously, there are at least two such points (since $C$ is infinite).
\end{proof}
Note that in the proof of Theorem \ref{gaga}, we only used the fact that the set $C$ is infinite. So, a similar result holds for any infinite set $C$.
\begin{corollary}\label{gaga1}
Let $X$ be a compact metric space and let  $G$ be a closed relation on $X$.  If no periodic point is generated by $G$, then no Cantor set is finitely generated by $G$. 
\end{corollary}
\begin{proof}
The corollary follows directly from Theorem \ref{gaga}.
\end{proof}
\begin{definition}
Let $X$ be a compact metric space and let  $G$ be a closed relation on $X$. We say that the set $G$ is \emph{ \color{blue}  irrationally embedded} or \emph{ \color{blue}  i-embedded into $X\times X$}, if 
\begin{enumerate}
\item  the entropy of $G$ is positive,
\item  no periodic point is generated by $G$.
\end{enumerate}
\end{definition}
\begin{definition}
Let $X$ be a compact metric space and let  $G$ be a closed relation on $X$. We say that the set $G$ is \emph{ \color{blue}  almost  i-embedded into $X\times X$}, if 
\begin{enumerate}
\item  the entropy of $G$ is positive,
\item  exactly one periodic point is generated by $G$.
\end{enumerate}
\end{definition}
\section{Topological conjugation for closed relations on compact metric spaces}\label{s3}
In this section, we generalize the notion of topological conjugacy for continuous functions to the topological conjugacy for closed relations $G$ and $H$ on $X$ and $Y$, respectively, and show that periodic points, finitely generated Cantor sets, Mahavier products and the entropy of closed relations are  preserved by a topological conjugation. Topological conjugacy for closed relations was introduced by E.~Akin  in  his book \cite{A}, where more properties of such a conjugacy may be found. Some of the results that follow are well-known to the researchers in the topic; since their proofs are short, we give them anyway.
\begin{definition}
Let $X$ and $Y$ be metric spaces, and let $f:X\rightarrow X$ and $g:Y\rightarrow Y$ be any functions.  If there is a homeomorphism $\varphi:X\rightarrow Y$ such that 
$$
\varphi \circ f=g\circ \varphi,
$$
then we say that \emph{ \color{blue}  $f$ and $g$ are topological conjugates}. 
\end{definition}
The following is a well-known result that we generalize later in Theorem \ref{joj111}.
\begin{theorem}\label{entropyconjugate}
Let $X$ and $Y$ be  compact metric spaces and let $f:X\rightarrow X$ and $g:Y\rightarrow Y$ be continuous functions that are  topological conjugates.  Then 
$$
h(f)=h(g)
$$
and 
$$
\varprojlim(X,f) \textup{ is homeomorphic to } \varprojlim(Y,g).
$$
\end{theorem}
\begin{proof}
The proof of the first part of the statement can be found in \cite[Theorem 7.2, page 167]{Walters book} and the proof of the second part of the statement in \cite[Theorem 33, page 23]{IM1}.
\end{proof}
\begin{observation}
Let $X$ and $Y$ be metric spaces, and let $f:X\rightarrow X$ and $g:Y\rightarrow Y$ be any functions.  The following statements are equivalent.
\begin{enumerate}
\item The functions $f$ and $g$ are topological conjugates. 
\item There is a homeomorphism $\varphi:X\rightarrow Y$  such that for each $(x,y)\in X\times X$, the following holds:
$$
(x,y)\in \Gamma(f) \Longleftrightarrow (\varphi(x),\varphi(y))\in \Gamma(g).
$$
\end{enumerate}
\end{observation}
Next, we generalize the notion of topological conjugacy of continuous functions to the  topological conjugacy of closed relations.
\begin{definition}
Let $X$ and $Y$ be any compact metric spaces, let $G$ be a closed relation on $X$ and let $H$ be a closed relation on $Y$. We say that \emph{ \color{blue}  $G$ and $H$ are topological conjugates} if there is a homeomorphism $\varphi:X\rightarrow Y$ such that for each $(x,y)\in X\times X$, the following holds
 $$
 (x,y)\in G  \Longleftrightarrow (\varphi(x), \varphi(y))\in H.
 $$
\end{definition}
The following theorem shows that the concept of conjugacy of upper semicontinuous functions, defined by Ingram and Mahavier in \cite{IM2}, is the same as the one that is introduced in the above definition for $G$ and $H$ being the graphs of upper semicontinuous set-valued functions.
\begin{theorem}
Let $X$ and $Y$ be compact metric spaces, and let $F:X\multimap X$ and $G:Y\rightarrow Y$ be upper semicontinuous set-valued functions. Also, let $\varphi:X\rightarrow Y$ be a homeomorphism. The following statements are equivalent.
\begin{enumerate}
\item For each $x\in X$, $G(\varphi(x))=\varphi(F(x))$.
\item For each $(x,y)\in X\times X$, 
$$
(x,y)\in \Gamma(F) \Longleftrightarrow (\varphi(x),\varphi(y))\in \Gamma(G).
$$
\end{enumerate}
\end{theorem}
\begin{proof}
Suppose that for each $x\in X$, $G(\varphi(x))=\varphi(F(x))$ and let $(x,y)\in X\times X$.
Suppose that $(x,y)\in \Gamma(F)$. Then $y\in F(x)$ and, therefore,  $\varphi(y)\in \varphi(F(x))$. Since $\varphi(F(x))=G(\varphi(x))$, it follows that $\varphi(y)\in G(\varphi(x))$. Therefore, $(\varphi(x),\varphi(y))\in \Gamma(G)$. Next, suppose that $(\varphi(x),\varphi(y))\in \Gamma(G)$.  It follows that $\varphi(y)\in G(\varphi(x))$ and, therefore, $y\in \varphi^{-1}(G(\varphi(x)))$.  Since $G(\varphi(x))=\varphi(F(x))$, it follows that $y\in \varphi^{-1}(\varphi(F(x)))$. Therefore, $y\in F(x)$ and $(x,y)\in \Gamma(F)$ follows. 

Next, suppose that for each $(x,y)\in X\times X$, 
$$
(x,y)\in \Gamma(F) \Longleftrightarrow (\varphi(x),\varphi(y))\in \Gamma(G)
$$
and let $x\in X$.  To show that $G(\varphi(x))\subseteq \varphi(F(x))$, let $y\in G(\varphi(x))$. Then $(\varphi(x),y)\in \Gamma(G)$.  It follows that $(x,\varphi^{-1}(y))\in \Gamma(F)$, and, therefore, $\varphi^{-1}(y)\in F(x)$.  It follows that $y\in \varphi(F(x))$. Finally, to show that $G(\varphi(x))\supseteq \varphi(F(x))$, let $y\in \varphi(F(x))$.  Then $\varphi^{-1}(y)\in F(x)$ and $(x,\varphi^{-1}(y))\in \Gamma(F)$ follows.  Hence, $(\varphi(x),y)\in \Gamma(G)$ and it follows that $y\in G(\varphi(x))$.  This completes the proof.    
\end{proof}
Next, we give results showing  that topological conjugations preserve finitely generated Cantor sets and periodic points. 
\begin{theorem}\label{joj2}
Let $X$ and $Y$ be any compact metric spaces, $G$ a closed relation on $X$ and $H$ a closed relation on $Y$. If $G$ and $H$ are topological conjugates, then the following statements are equivalent.
\begin{enumerate}
\item  No Cantor set is finitely generated by $G$.
\item No Cantor set is finitely generated by $H$.
\end{enumerate}
\end{theorem}
\begin{proof}
Let $\varphi:X\rightarrow Y$ be a homeomorphism such that for each $(x,y)\in X\times X$, 
 $$
 (x,y)\in G  \Longleftrightarrow (\varphi(x), \varphi(y))\in H.
 $$
Suppose that $C\subseteq \star_{n=1}^{\infty}G^{-1}$ is a Cantor set and that 
$$
F_1=\{(x_1,y_1),(x_2,y_2),(x_3,y_3),\ldots, (x_n,y_n)\}\subseteq G^{-1}
$$
 such that 
$$
C=\star_{n=1}^{\infty}F_1.
$$
Let
$$
F_2=\{(\varphi(x_1),\varphi(y_1)),(\varphi(x_2),\varphi(y_2)),(\varphi(x_3),\varphi(y_3)),\ldots, (\varphi(x_n),\varphi(y_n))\}.
$$
Then  $F_2\subseteq H^{-1}$.
The function
$$
\Phi : \star_{n=1}^{\infty}F_1\rightarrow \star_{n=1}^{\infty}F_2, 
$$
defined by 
$$
\Phi(x_1,x_2,x_3,\ldots)=(\varphi(x_1),\varphi(x_2),\varphi(x_3),\ldots)
$$
for each $(x_1,x_2,x_3,\ldots)\in \star_{n=1}^{\infty}F_1$ is  a homeomorphism. Therefore,   $\star_{n=1}^{\infty}F_2$ is a Cantor set that is finitely generated by $H$.
 
 Next, suppose that $D\subseteq \star_{n=1}^{\infty}H^{-1}$ is a Cantor set and that 
$$
H_1=\{(x_1,y_1),(x_2,y_2),(x_3,y_3),\ldots, (x_n,y_n)\}\subseteq H^{-1}
$$
 such that 
$$
D=\star_{n=1}^{\infty}H_1.
$$
Let
$$
H_2=\{(\varphi^{-1}(x_1),\varphi^{-1}(y_1)),(\varphi^{-1}(x_2),\varphi^{-1}(y_2)),(\varphi^{-1}(x_3),\varphi^{-1}(y_3)),\ldots, (\varphi^{-1}(x_n),\varphi^{-1}(y_n))\}.
$$
Then  $H_2\subseteq G^{-1}$.
The function
$$
\Phi : \star_{n=1}^{\infty}H_1\rightarrow \star_{n=1}^{\infty}H_2, 
$$
defined by 
$$
\Phi(x_1,x_2,x_3,\ldots)=(\varphi^{-1}(x_1),\varphi^{-1}(x_2),\varphi^{-1}(x_3),\ldots)
$$
for each $(x_1,x_2,x_3,\ldots)\in \star_{n=1}^{\infty}H_1$, is  a homeomorphism.  Therefore,  $\star_{n=1}^{\infty}H_2$ is a Cantor set that is finitely generated by $G$.
\end{proof}
\begin{lemma}\label{joj0}
Let $X$ and $Y$ be any compact metric spaces, let $G$ be a closed relation on $X$ and let $H$ be a closed relation on $Y$,  and let $(x_1,x_2,x_3,\ldots)\in \star_{i=1}^{\infty}G^{-1}$. If $\varphi:X\rightarrow Y$ is a homeomorphism such that for each $(x,y)\in X\times X$, 
 $$
 (x,y)\in G  \Longleftrightarrow (\varphi(x), \varphi(y))\in H,
 $$
  then the following statements are equivalent.
\begin{enumerate}
\item\label{1}  $(x_1,x_2,x_3,\ldots)$ is a periodic point with period $m$ that is generated by $G$.
\item \label{2} $(\varphi(x_1),\varphi(x_2),\varphi(x_3),\ldots)$ is a periodic point with period $m$  that is generated by $H$.
\end{enumerate}
\end{lemma}
\begin{proof}
Let $\sigma_G:\star_{n=1}^{\infty}G^{-1}\rightarrow \star_{n=1}^{\infty}G^{-1}$ be the shift function on $\star_{n=1}^{\infty}G^{-1}$ and let $\sigma_H:\star_{n=1}^{\infty}H^{-1}\rightarrow \star_{n=1}^{\infty}H^{-1}$ be the shift function on $\star_{n=1}^{\infty}H^{-1}$. 
Suppose that $\mathbf x=(x_1,x_2,x_3,\ldots)\in \star_{n=1}^{\infty}G^{-1}$ is a periodic point that is generated by $G$.  Let $m$ be a positive integer such that 
$\sigma_G^{m}(\mathbf x)=\mathbf x$.  Let 
$$
\mathbf y=(\varphi(x_1),\varphi(x_2),\varphi(x_3),\ldots).
$$
Obviously, $\mathbf y\in \star_{n=1}^{\infty}H^{-1}$. We show that $\sigma_H^{m}(\mathbf y)=\mathbf y$. Since $\sigma_G^{m}(\mathbf x)=\mathbf x$, it follows that for each positive integer $i$,
$$
x_{i}=x_{i+m}.
$$
Therefore, for each positive integer $i$,
$$
\varphi(x_i)=\varphi(x_{i+m})
$$
and it follows that $\sigma_H^{m}(\mathbf y)=\mathbf y$.  This completes the proof that \ref{1} implies \ref{2}.

To show that \ref{2} implies \ref{1}, suppose that $\mathbf y=(\varphi(x_1),\varphi(x_2),\varphi(x_3),\ldots)\in \star_{n=1}^{\infty}H^{-1}$ is a periodic point that is generated by $H$.  Let $m$ be a positive integer such that 
$\sigma_H^{m}(\mathbf y)=\mathbf y$.  Let $\mathbf x=(x_1,x_2,x_3,\ldots)$. 
 Since $\sigma_H^{m}(\mathbf y)=\mathbf y$, it follows that for each positive integer $i$,
$$
\varphi(x_{i})=\varphi(x_{i+m}).
$$
Therefore, for each positive integer $i$,
$$
x_i=x_{i+m}
$$
and it follows that $\sigma_G^{m}(\mathbf x)=\mathbf x$.  
\end{proof}
\begin{theorem}\label{joj1}
Let $X$ and $Y$ be any compact metric spaces, $G$ a closed relation on $X$ and $H$ a closed relation on $Y$. If $G$ and $H$ are topological conjugates, then the following statements are equivalent.
\begin{enumerate}
\item  No periodic point is generated by $G$.
\item No periodic point is generated by $H$.
\end{enumerate}
\end{theorem}
\begin{proof}
The result follows directly from Lemma \ref{joj0}.
\end{proof}

\begin{theorem}\label{joj111}
Let $X$ and $Y$ be any compact metric spaces, $G$ a closed relation on $X$ and $H$ a closed relation on $Y$. If $G$ and $H$ are topological conjugates, then 
$$
\ent(G)=\ent(H) 
$$
and
$$
\star_{i=1}^{\infty}G^{-1} \textup{ is homeomorphic to } \star_{i=1}^{\infty}H^{-1}.
$$
\end{theorem}
\begin{proof}
Let $\varphi:X\rightarrow Y$ be a homeomorphism such that for each $(x,y)\in X\times X$, 
 $$
 (x,y)\in G  \Longleftrightarrow (\varphi(x), \varphi(y))\in H.
 $$
Note that for any open cover  $\alpha$ for $X$,  $\varphi(\alpha)=\{\varphi(U) \ | \ U\in \alpha\}$ is an open cover for $Y$, and that for any open cover  $\beta$ for $Y$,  $\varphi^{-1}(\beta)=\{\varphi^{-1}(V) \ | \ V\in \beta\}$ is an open cover for $X$. 
Let $\alpha$ be any open cover for $X$, let $N(\star_{i=1}^mG^{-1},\alpha^{m+1})=n$ and let $\mathcal U=\{U_1,U_2,U_3,\ldots,U_n\}$ be a subcover of $\alpha^{m+1}$ such that $\star_{i=1}^mG^{-1}\subseteq U_1\cup U_2\cup U_3\cup \ldots\cup U_n$.   Then $\star_{i=1}^mH^{-1}\subseteq \varphi(U_1)\cup \varphi(U_2)\cup \varphi(U_3)\cup \ldots\cup \varphi(U_n)$ and
$$
N(\star_{i=1}^{m}H^{-1},\varphi(\alpha)^{m+1})\leq  N(\star_{i=1}^{m}G^{-1},\alpha^{m+1})
$$
follows. Next, let $\beta$ be any open cover for $Y$, let $N(\star_{i=1}^mH^{-1},\beta^{m+1})=n$ and let $\mathcal V=\{V_1,V_2,V_3,\ldots,V_n\}$ be a subcover of $\beta^{m+1}$ such that $\star_{i=1}^mH^{-1}\subseteq V_1\cup V_2\cup V_3\cup \ldots\cup V_n$.   Then $\star_{i=1}^mG^{-1}\subseteq \varphi^{-1}(V_1)\cup \varphi^{-1}(V_2)\cup \varphi^{-1}(V_3)\cup \ldots\cup \varphi^{-1}(V_n)$.  Therefore,
$$
N(\star_{i=1}^{m}G^{-1},\varphi^{-1}(\beta)^{m+1})\leq N(\star_{i=1}^{m}H^{-1},\beta^{m+1}).
$$
It follows that for each open cover  $\alpha$ for $X$,
$$
\ent(G,\alpha) =\lim_{m\to \infty}\frac{N(\star_{i=1}^{m}G^{-1},\alpha^{m+1})}{m}\geq \lim_{m\to \infty}\frac{N(\star_{i=1}^{m}H^{-1},\varphi(\alpha)^{m+1})}{m}=\ent(H,\varphi(\alpha))
$$
and that for each open cover  $\beta$ for $Y$,
$$
\ent(H,\beta) =\lim_{m\to \infty}\frac{N(\star_{i=1}^{m}H^{-1},\beta^{m+1})}{m}\geq \lim_{m\to \infty}\frac{N(\star_{i=1}^{m}G^{-1},\varphi^{-1}(\beta)^{m+1})}{m}=\ent(G,\varphi^{-1}(\beta)).
$$
Therefore, 
 $$
 \ent(G)=\sup\{\ent(G,\alpha) \ | \ \alpha \textup{ is an open cover for } X\} 
 $$
 $$
 =\sup\{\ent(H,\beta) \ | \ \beta \textup{ is an open cover for } Y\} =\ent(H).
 $$
 Also, note that the function 
 $$
\psi:\star_{i=1}^{\infty}G^{-1} \rightarrow \star_{i=1}^{\infty}H^{-1},
$$
 defined by 
 $$
 \psi(x_1,x_2,x_3,\ldots)=(\varphi(x_1), \varphi(x_2), \varphi(x_3), \ldots)
 $$
 for any $(x_1,x_2,x_3,\ldots)\in \star_{i=1}^{\infty}G^{-1}$, is a homeomorphism from $\star_{i=1}^{\infty}G^{-1}$ to $\star_{i=1}^{\infty}H^{-1}$. This completes our proof.
\end{proof}
\section{A sufficient condition for  non-zero entropy of closed relations on $[0,1]$}\label{s4}
\begin{figure}[h!]
	\centering
		\includegraphics[width=30em]{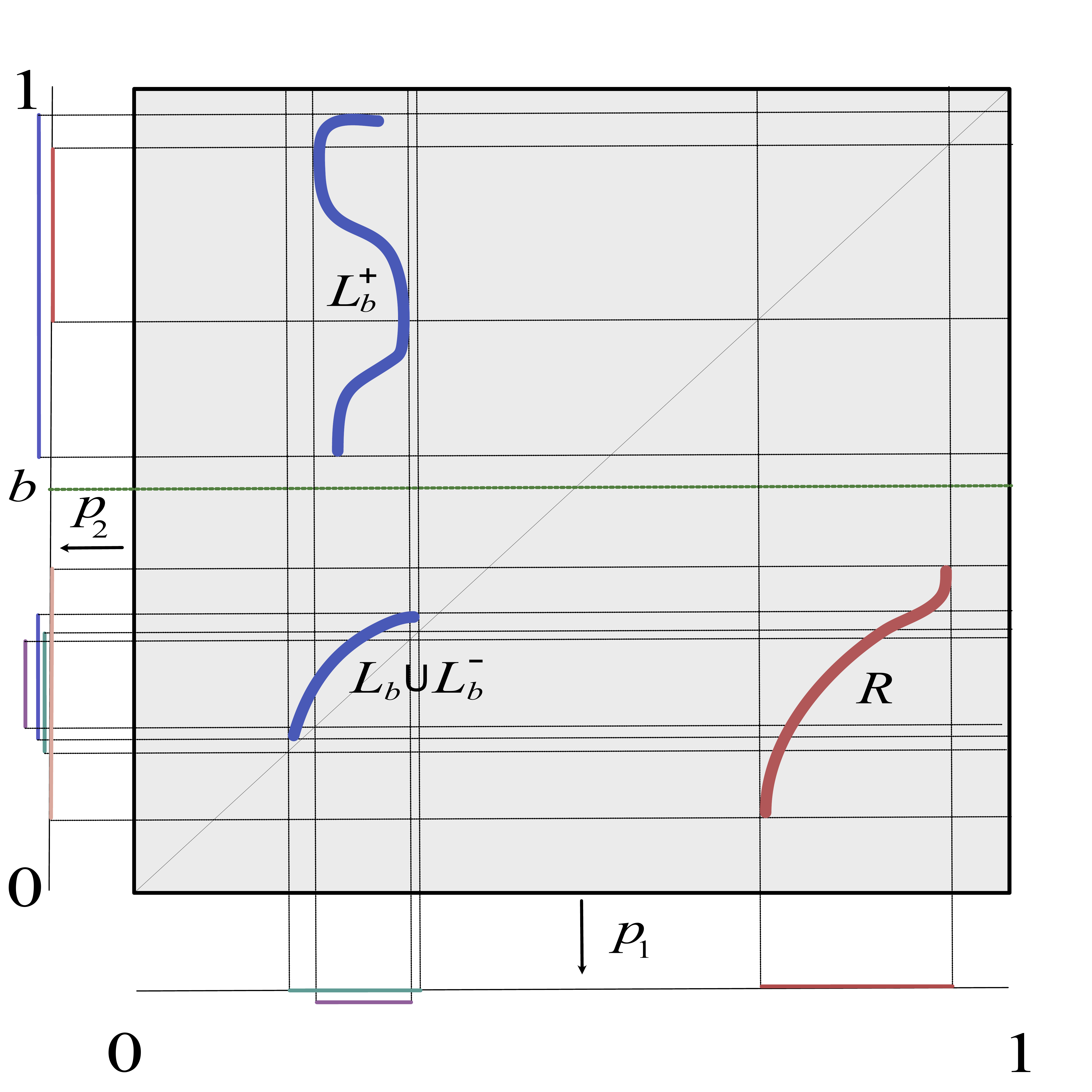}
	\caption{The sets $L$ and $R$}
	\label{figure12321}
\end{figure}

In the rest of the paper,  $I$ will always denote the closed interval $[0,1]$. We also use 
$$
p_1:I\times I\rightarrow I
$$
and 
$$
p_2:I\times I\rightarrow I
$$
to denote the standard projections defined by
$$
p_1(s,t)=s  \textup{ and } p_2(s,t)=t
$$
for all $(s,t)\in I\times I$.
\begin{definition}
We define
$$
\Delta^{+}=\{(x,y)\in I\times I \ | \ y>x\},
$$
$$
\Delta^{-}=\{(x,y)\in I\times I \ | \ y<x\}
$$
and
$$
\Delta=\{(x,y)\in I\times I \ | \ y=x\}.
$$
\end{definition}
\begin{definition}
Let $A\subseteq I\times I$ and $b\in (0,1)$. Then we define
$$
A^{+}_{b}=\{(x,y)\in A \ | \ y>b\},
$$
$$
A^{-}_{b}=\{(x,y)\in A \ | \ y<b\},
$$
and
$$
A_{b}=\{(x,y)\in A \ | \ y=b\}.
$$
\end{definition}
\begin{definition}
Let $L$ and $R$ be closed subsets of $I\times I$ and let $b\in (0,1)$. We say that the sets \emph{ \color{blue}  $L$ and $R$ are well-aligned by $b$},  if 
\begin{enumerate}
\item $L_{b}^{+}\neq \emptyset$, $L_{b}^{-}\cup L_b\neq \emptyset$ and $R=R_b\cup R_b^{-}$,
\item $L_{b}^{+}\cup L_b\subseteq \Delta^{+}$, $L_{b}^{-}\subseteq \Delta^{+}\cup \Delta$ and $R\subseteq \Delta^{-}$,
\item $p_2(L_b^{-}\cup L_b)\cup p_1(L_b^{-}\cup L_b)\subseteq p_2(R)$,
\item $p_1(L_b^{+})\cup p_1(R)\subseteq p_2(L)$,
\end{enumerate}
see Figure \ref{figure12321}.
\end{definition}
\begin{definition}
Let $L$ and $R$ be  closed subsets of $I\times I$.  We say that the sets \emph{ \color{blue}  $L$ and $R$ are well-aligned,} if there is a point  $b\in (0,1)$ such that the sets $L$ and $R$ are well-aligned by $b$. 
\end{definition}

\begin{definition}
Let $G$ be a non-empty closed subset of $I\times I$.
Then we define the functions $r_{G},\ell_{G}:p_2(G)\rightarrow p_1(G)$ by
$$
\ell_{G}(t)=\min\{p_1(s,t) \ | \ (s,t)\in G\}
$$
and
$$
r_{G}(t)=\max\{p_1(s,t) \ | \ (s,t)\in G\}
$$
for each $t\in p_2(G)$.
\end{definition}
\begin{lemma}\label{mini}
Let $b\in (0,1)$ and let $L$ and $R$ be closed subsets of $I\times I$ that are well-aligned by $b$. 
Then the following holds.
\begin{enumerate}
\item For each $t\in p_2(L_b^{+})$ there is a positive integer $k_t$ such that
$$
r_{L}^{k_t}(t)\leq b<r_{L}^{k_t-1}(t).
$$
\item For each $t\in p_2(L_b^{+})$ let $k_t$ be the  positive integer such that 
$$
r_{L}^{k_t}(t)\leq b<r_{L}^{k_t-1}(t). 
$$
Then there is a positive integer $k$ such that for each $t\in p_2(L_b^{+})$, 
$$
k_t\in \{1,2,3,\ldots , k\}.
$$
\end{enumerate}
\end{lemma}
\begin{proof}
Let $d$ be the Euclidean metric on $I\times I$. It follows from $L_{b}^{+}\cup L_b\subseteq \Delta^{+}$ that $d(L_{b}^{+}\cup L_b,\Delta)>0$.  Therefore, there is an $a<1$ such that for the function $f:I\rightarrow I$, for each $t\in I$ defined by 
$$
f(t)=a\cdot t,
$$
 the following holds: 
$$
r_{L}(t)\leq f(t)<t
$$
 for each $t\in p_2(L_{b}^{+}\cup L_b)$.  Fix such an $a$ and such a linear function $f$. 

Let $k$ be the positive integer such that $f^{k}(1) \leq  b<f^{k-1}(1)$.  Let $t\in p_2(L_b^{+})$. We show that there is an integer $k_t\in \{1,2,3,\ldots,k\}$ such that  
$$
r_{L}^{k_t}(t)\leq b<r_{L}^{k_t-1}(t).
$$ 
Note that $(r_{L}(t),t)\in L_b^{+}$ (since $t\in p_2(L_b^{+})$), therefore,  $r_{L}(t)\in p_1(L_b^{+})$.  Since $p_1(L_b^{+})\subseteq p_2(L)$, it follows that  $r_{L}(t)\in p_2(L)$. If $r_{L}(t)\in p_2(L_{b}^{-}\cup L_b)$, we are done. Otherwise, $r_{L}(t)\in p_2(L_{b}^{+})$.  Then $(r_{L}^{2}(t),r_{L}(t))\in L_b^{+}$, therefore,  $r_{L}^{2}(t)\in p_1(L_b^{+})$.  and $r_{L}^{2}(t)\in p_2(L)$ follows. If $r_{L}^{2}(t)\in p_2(L_{b}^{-}\cup L_b)$, we are done. Otherwise, $r_{L}^{2}(t)\in p_2(L_{b}^{+})$.  A similar argument shows that also in this case  $r_{L}^{3}(t)\in p_2(L)$. If $r_{L}^{3}(t)\in p_2(L_{b}^{-}\cup L_b)$, we are done. Otherwise, $r_{L}^{3}(t)\in p_2(L_{b}^{+})$.  We continue inductively. Since $f^{k}(t) \leq  b$ and $r_{L}(t)\leq f(t)$, it follows that there is an integer $k_t\in \{1,2,3,\ldots,k\}$ such that  
$$
r_{L}^{k_t}(t)\leq b<r_{L}^{k_t-1}(t).
$$ 
\end{proof}
\begin{definition}
Let $b\in (0,1)$ and let $L$ and $R$ be closed subsets of $I\times I$ that are well-aligned by $b$,  and let $t\in p_2(L)$ be any point.  If $t\in p_2(L_b^{+})$, then we define $\Psi^{L,R}_{b}(t)$ to be the positive integer $k$ such that  
$$
r_{L}^{k}(t)\leq b<r_{L}^{k-1}(t).
$$
If $t\in p_2(L_b^{-}\cup L_b)$, then we define $\psi^{L,R}_{b}(t)=0$.
We also define
$$
\psi^{L,R}_{b}=\max \{\psi^{L,R}_{b}(t) \ | \ t\in p_2(L)\}.
$$
\end{definition}
\begin{definition}
Let $X$ be a set, let $(i_k)$ be a sequence of positive integers, and let $T_k=(x_{k,1},x_{k,2},x_{k,3},\ldots ,x_{k,i_k})\in X^{i_k}$ for each positive integer $k$. We define
$$
T_1\oplus T_2\oplus T_3\oplus\ldots \oplus T_n =\oplus_{k=1}^{n}T_k=
$$
$$
(x_{1,1},x_{1,2},x_{1,3},\ldots ,x_{1,i_1},x_{2,1},x_{2,2},x_{2,3},\ldots ,x_{2,i_2},\ldots, x_{n,1},x_{n,2},x_{n,3},\ldots ,x_{n,i_n})
$$
and
$$
T_1\oplus T_2\oplus T_3\oplus\ldots =\oplus_{k=1}^{\infty}T_k=(x_{1,1},x_{1,2},x_{1,3},\ldots ,x_{1,i_1},x_{2,1},x_{2,2},x_{2,3},\ldots ,x_{2,i_2},\ldots)
$$
\end{definition}

\begin{theorem}\label{MAIN}
Let $G$ be a closed relation on $I$.  If there are well-aligned  closed subsets of $G$ or if there are well-aligned  closed subsets of $G^{-1}$, then $\ent(G)\neq 0$. 
\end{theorem}
\begin{proof}
Let $b\in (0,1)$, and let  $L$  and $R$ be closed subsets of $G$ such that
\begin{enumerate}
\item $L_{b}^{+}\neq \emptyset$, $L_{b}^{-}\cup L_b\neq \emptyset$ and $R=R_b\cup R_b^{-}$,
\item $L_{b}^{+}\cup L_b\subseteq \Delta^{+}$, $L_{b}^{-}\subseteq \Delta^{+}\cup \Delta$ and $R\subseteq \Delta^{-}$,
\item $p_2(L_b^{-}\cup L_b)\cup p_1(L_b^{-}\cup L_b)\subseteq p_2(R)$,
\item $p_1(L_b^{+})\cup p_1(R)\subseteq p_2(L)$.
\end{enumerate}
Let 
 $$
 \varepsilon=\inf\{\ell_R(t)-r_L(t) \ | \ t\in p_2(L)\cap p_2(R)\}
 $$
 and let $\alpha_0$ be any open cover for $I$ such that $\diam(A)<\frac{\varepsilon}{2}$ for any $A\in \alpha_0$.  
 
 For each $t\in p_2(L^{+}_{b})$, let 
 $$
 T(t)=\left(t,r_{L}(t),r_{L}^{2}(t),r_{L}^{3}(t),\ldots,r_{L}^{\psi^{L,R}_{b}(t)-1}(t), r_{L}^{\psi^{L,R}_{b}(t)}(t)\right).
 $$
Let $t\in p_2(L^{+}_{b})$ be any point and let $i_1=\psi^{L,R}_{b}(t)$, ${\color{red}t_0=r_{L}^{i_1+1}(t)}$   and ${\color{red}t_1=\ell_{R}(r_{L}^{i_1}(t))}$. 
Then, obviously,  $T(t)\oplus {\color{red}t_0}$ and $T(t)\oplus {\color{red}t_1}$ are both elements of $\star_{i=1}^{i_1+1}G^{-1}$.
It follows that
$$ 
{\color{red}t_1}-{\color{red}t_0}\geq \inf\{\ell_R(t)-r_L(t) \ | \ t\in p_2(L)\cap p_2(R)\}=\varepsilon
$$
 holds. Therefore,  $\{{\color{red}t_0},{\color{red}t_1}\}\not \subseteq A$ since $\diam(A)<\frac{\varepsilon}{2}$ for any $A\in \alpha_0$.   
It follows that for any $U\in \alpha_0^{i_1+2}$,
$$
 \{T(t)\oplus {\color{red}t_0},T(t)\oplus {\color{red}t_1}\}\not \subseteq U
 $$
and, therefore, 
$$
N(\star_{i=1}^{i_1+1}G^{-1},\alpha_0^{i_1+2})\geq 2.
$$

It follows from $r_{L}^{i_1}(t)\in p_2(L_b^{-}\cup L_b)$ that $ t_0\in p_1(L_b^{-}\cup L_b)$.  Since $p_1(L_b^{-}\cup L_b)\subseteq p_2(R)$, it follows that $ t_0\in p_2(R)$.  Therefore, $(\ell_{R}(t_0),t_0)\in R$ and $\ell_{R}(t_0)\in p_1(R)$ follows.  Since $p_1(R)\subseteq p_2(L)$, it follows that $\ell_{R}(t_0)\in p_2(L)$. Let 
$$
\psi_{0}=\psi^{L,R}_b(\ell_{R}(t_0)).
$$
Note that $\psi_{0}\in \{0,1,2,3,\ldots, \psi^{L,R}_b\}$. Also, let $T_0=t_0\oplus T(\ell_{R}(t_0))$, ${\color{blue} t_{00}=r_{L}^{\psi_{0}+1}(\ell_{R}(t_0))}$ and ${\color{blue} t_{01}=\ell_R(r_{L}^{\psi_{0}}(\ell_{R}(t_0)))}$.
Then
$$
T(t)\oplus T_0\oplus {\color{blue} t_{00}}
$$ 
and
$$
T(t)\oplus T_0\oplus {\color{blue} t_{01}}
$$
are both elements of $\star_{i=1}^{i_1+\psi_{0}+3}G^{-1}$.  It follows that
$$ 
{\color{blue}t_{01}}-{\color{blue}t_{00}}\geq \inf\{\ell_R(t)-r_L(t) \ | \ t\in p_2(L)\cap p_2(R)\}=\varepsilon
$$
 holds.  Therefore,  $\{t_{00},t_{01}\}\not \subseteq A$ since $\diam(A)<\frac{\varepsilon}{2}$ for any $A\in \alpha_0$.   
It follows that for any $U\in \alpha_0^{i_1+\psi_{0}+4}$,
$$
 \{T(t)\oplus T_0\oplus {\color{blue} t_{00}},T(t)\oplus T_0\oplus {\color{blue} t_{01}}\}\not \subseteq U.
 $$ 
 Next,  since $t_1\in p_1(R)$ and $p_1(R)\subseteq p_2(L)$, it follows that $ t_1\in p_2(L)$.  Let 
$$
\psi_{1}=\psi^{L,R}_b(t_1).
$$
Note that $\psi_{1}\in \{0,1,2,3,\ldots, \psi^{L,R}_b\}$. Also, let $T_1=T(t_1)$, ${\color{green} t_{10}=r_{L}^{\psi_{1}+1}(t_1)}$ and ${\color{green} t_{11}=\ell_R(r_{L}^{\psi_{1}}(t_1))}$.
Then
$$
T(t)\oplus T_1\oplus {\color{green} t_{10}}
$$ 
and
$$
T(t)\oplus T_1\oplus {\color{green} t_{11}}
$$
are both elements of $\star_{i=1}^{i_1+\psi_{1}+2}G^{-1}$.  It follows that
$$ 
{\color{green}t_{11}}-{\color{green}t_{10}}\geq \inf\{\ell_R(t)-r_L(t) \ | \ t\in p_2(L)\cap p_2(R)\}=\varepsilon
$$
 holds.  Therefore,  $\{t_{10},t_{11}\}\not \subseteq A$ since $\diam(A)<\frac{\varepsilon}{2}$ for any $A\in \alpha_0$.   
It follows that for any $U\in \alpha_0^{i_1+\psi_{1}+3}$,
$$
 \{T(t)\oplus T_1\oplus {\color{green} t_{10}},T(t)\oplus T_1\oplus {\color{green} t_{11}}\}\not \subseteq U.
 $$ 
 Let  
$$
i_2=\max\{\psi_{0},\psi_{1}\}.
$$
Let $P_{00},P_{01},P_{10},P_{11}$ be any points in $\star_{i=1}^{i_1+i_2+3}G^{-1}$ such that 
$$
p_{[1,i_1+\psi_{0}+4]}(P_{00})=T(t)\oplus T_0\oplus {\color{blue} t_{00}},
$$
$$
p_{[1,i_1+\psi_{0}+4]}(P_{01})=T(t)\oplus T_0\oplus {\color{blue} t_{01}},
$$
$$
p_{[1,i_1+\psi_{1}+3]}(P_{10})=T(t)\oplus T_1\oplus {\color{green} t_{10}},
$$
and
$$
p_{[1,i_1+\psi_{1}+3]}(P_{11})=T(t)\oplus T_1\oplus {\color{green} t_{11}},
$$
Note that for each $U\in \alpha_0^{i_1+i_2+4}$ and for any two points $\mathbf x,\mathbf y\in \{P_{00},P_{01},P_{10},P_{11}\}$ it holds that if $\mathbf x\neq \mathbf y$, then 
$$
\{\mathbf x, \mathbf y\}\not \subseteq  U.
$$
 It follows that 
$$
N(\star_{i=1}^{i_1+i_2+3}G^{-1},\alpha_0^{i_1+i_2+4})\geq 2^2.
$$
We continue inductively. Let $n$ be a positive integer and suppose that for each positive integer $k\leq n$ and  for each $(z_1,z_2,z_3,\ldots,z_k)\in \{0,1\}^{k}$, we have already constructed 
\begin{enumerate}
\item the positive integer 
$$
\psi_{z_1z_2z_3\ldots z_k}\in  \{0,1,2,3,\ldots, \psi^{L,R}_b\},
$$
\item the point $T_{z_1z_2z_3\ldots z_k}$ such that $T_{z_1z_2z_3\ldots z_k}\in X^{\psi_{z_1z_2z_3\ldots z_k}+2}$, if $z_k=0$, and $T_{z_1z_2z_3\ldots z_k}\in X^{\psi_{z_1z_2z_3\ldots z_k}+1}$, if $z_k=1$,
\item the points $t_{z_1z_2z_3\ldots z_k0}\in p_1(L)$ and $t_{z_1z_2z_3\ldots z_k1}\in p_1(R)$ such that 
$$ 
t_{z_1z_2z_3\ldots z_k1}-t_{z_1z_2z_3\ldots z_k0}\geq \inf\{\ell_R(t)-r_L(t) \ | \ t\in p_2(L)\cap p_2(R)\}=\varepsilon,
$$
 and
\item the points $P_{z_1z_2z_3\ldots z_k0}$ and $P_{z_1z_2z_3\ldots z_k1}$ in $\star_{i=1}^{i_1+i_2+i_3+\ldots +i_k+i_{k+1}+2k+1}G^{-1}$,
\end{enumerate}
where $i_{k+1}=\max\{\psi_{z_1z_2z_3\ldots z_k} \ | \ (z_1,z_2,z_3,\ldots,z_k)\in \{0,1\}^{k}\}$, such that for each $z\in \{0,1\}$,
$$
p_{[1,i_1+\psi_{z_1}+\psi_{z_1z_2}+\psi_{z_1z_2z_3}+\ldots+\psi_{z_1z_2z_3\ldots z_k}+2k+2-w]}(P_{z_1z_2z_3\ldots z_kz})
$$
$$
=T(t)\oplus T_{z_1} \oplus T_{z_1z_2}\oplus T_{z_1z_2z_3}\oplus \ldots \oplus T_{z_1z_2z_3\ldots z_k}\oplus  t_{z_1z_2z_3\ldots z_kz},
$$
where $w$ is the number of occurrences of a digit $1$ in the sequence $z_1z_2z_3\ldots z_k$; meaning that 
$$
N(\star_{i=1}^{i_1+i_2+i_3+\ldots +i_k+i_{k+1}+2k+1}G^{-1},\alpha_0^{i_1+i_2+i_3+\ldots +i_k+i_{k+1}+2k+2})\geq 2^{k+1}
$$
for each $k\leq n$.

Next, let $(z_1,z_2,z_3,\ldots,z_n,z_{n+1})\in \{0,1\}^{n+1}$. Then we define
\begin{enumerate}
\item the positive integer $\psi_{z_1z_2z_3\ldots z_nz_{n+1}}$ as follows:
$$
\psi_{z_1z_2z_3\ldots z_nz_{n+1}} = \begin{cases}
				\psi_{b}^{L,R}(\ell_R(t_{z_1z_2z_3\ldots z_nz_{n+1}}))\text{;} & z_{n+1}=0\\
				\psi_{b}^{L,R}(t_{z_1z_2z_3\ldots z_nz_{n+1}})\text{;} & z_{n+1}=1.
			\end{cases}
$$
Note that $\psi_{z_1z_2z_3\ldots z_nz_{n+1}}\in  \{0,1,2,3,\ldots, \psi^{L,R}_b\}$. 
\item the point $T_{z_1z_2z_3\ldots z_nz_{n+1}}$ as follows:
$$
T_{z_1z_2z_3\ldots z_nz_{n+1}} = \begin{cases}
				t_{z_1z_2z_3\ldots z_nz_{n+1}}\oplus T(\ell_R(t_{z_1z_2z_3\ldots z_nz_{n+1}}))\text{;} & z_{n+1}=0\\
				T(t_{z_1z_2z_3\ldots z_nz_{n+1}})\text{;} & z_{n+1}=1.
			\end{cases}
$$
\item the points $t_{z_1z_2z_3\ldots z_nz_{n+1}0}\in p_1(L)$ and $t_{z_1z_2z_3\ldots z_nz_{n+1}1}\in p_1(R)$ as follows;
$$
t_{z_1z_2z_3\ldots z_nz_{n+1}0}=r_L^{\psi_{z_1z_2z_3\ldots z_nz_{n+1}}+1}(\ell_R(t_{z_1z_2z_3\ldots z_nz_{n+1}}))
$$
and
$$
t_{z_1z_2z_3\ldots z_nz_{n+1}1}=\ell_R(r_L^{\psi_{z_1z_2z_3\ldots z_nz_{n+1}}}(\ell_R(t_{z_1z_2z_3\ldots z_nz_{n+1}}))),
$$
if $z_{n+1}=0$, and 
$$
t_{z_1z_2z_3\ldots z_nz_{n+1}0}=r_L^{\psi_{z_1z_2z_3\ldots z_nz_{n+1}}+1}(t_{z_1z_2z_3\ldots z_nz_{n+1}})
$$
and
$$
t_{z_1z_2z_3\ldots z_nz_{n+1}1}=\ell_R(r_L^{\psi_{z_1z_2z_3\ldots z_nz_{n+1}}}(t_{z_1z_2z_3\ldots z_nz_{n+1}})),
$$
if $z_{n+1}=1$.
Note that  
$$ 
t_{z_1z_2z_3\ldots z_nz_{n+1}1}-t_{z_1z_2z_3\ldots z_nz_{n+1}0}\geq \inf\{\ell_R(t)-r_L(t) \ | \ t\in p_2(L)\cap p_2(R)\}=\varepsilon,
$$
\end{enumerate}
Let $i_{n+2}=\max\{\psi_{z_1z_2z_3\ldots z_nz_{n+1}} \ | \ (z_1,z_2,z_3,\ldots,z_n,z_{n+1})\in \{0,1\}^{n+1}\}$.
For each 
$$
(z_1,z_2,z_3,\ldots,z_n,z_{n+1})\in \{0,1\}^{n+1},
$$
 let $P_{z_1z_2z_3\ldots z_nz_{n+1}0}$ and $P_{z_1z_2z_3\ldots z_nz_{n+1}1}$ be any points in $\star_{i=1}^{i_1+i_2+i_3+\ldots +i_n+i_{n+1}+i_{n+2}+2n+3}G^{-1}$,
such that for each $z\in \{0,1\}$,
$$
p_{[1,i_1+\psi_{z_1}+\psi_{z_1z_2}+\psi_{z_1z_2z_3}+\ldots+\psi_{z_1z_2z_3\ldots z_n}+\psi_{z_1z_2z_3\ldots z_nz_{n+1}}+2n+4-w]}(P_{z_1z_2z_3\ldots z_nz_{n+1}z})
$$
$$
=T(t)\oplus T_{z_1} \oplus T_{z_1z_2}\oplus T_{z_1z_2z_3}\oplus \ldots \oplus T_{z_1z_2z_3\ldots z_n}\oplus T_{z_1z_2z_3\ldots z_nz_{n+1}} \oplus  t_{z_1z_2z_3\ldots z_nz_{n+1}z},
$$
where $w$ is the number of occurrences of a digit $1$ in the sequence $z_1z_2z_3\ldots z_nz_{n+1}$. 
Note that for each $U\in \alpha_0^{i_1+i_2+i_3+\ldots +i_n+i_{n+1}+i_{n+2}+2n+4}$ and for any two points 
$$
\mathbf x,\mathbf y\in \{P_{z_1z_2z_3\ldots z_nz_{n+1}z_{n+2}} \ | \ (z_1,z_2,z_3,\ldots ,z_n,z_{n+1},z_{n+2})\in \{0,1\}^{n+2}\}
$$
 it holds that if $\mathbf x\neq \mathbf y$, then 
$$
\{\mathbf x,\mathbf y\}\not \subseteq  U.
$$
Therefore, 
$$
N(\star_{i=1}^{i_1+i_2+i_3+\ldots +i_n+i_{n+1}+i_{n+2}+2n+3}G^{-1},\alpha_0^{i_1+i_2+i_3+\ldots +i_n+i_{n+1}+i_{n+2}+2n+4})\geq 2^{n+2}.
$$
This completes the proof that for each positive integer $m\geq 2$,
$$
N(\star_{i=1}^{i_1+i_2+i_3+\ldots +i_m+2m-1}G^{-1},\alpha_0^{i_1+i_2+i_3+\ldots +i_m+2m})\geq 2^{m}.
$$
It follows that for each positive integer $m\geq 2$, 
$$
\frac{\log N(\star_{i=1}^{i_1+i_2+i_3+\ldots +i_m+2m-1}G^{-1},\alpha_0^{i_1+i_2+i_3+\ldots +i_m+2m})}{i_1+i_2+i_3+\ldots +i_m+2m-1}\geq  \frac{\log 2^m}{\psi^{L,R}_{b}\cdot m+2m}
$$
holds, since   $i_k\leq \psi^{L,R}_{b}$ for each $k\in \{1,2,3,\ldots,m\}$ and each positive integer $m$.  Therefore,
$$
\ent(G)=\sup\{\ent(G,\alpha) \ | \ \alpha \textup{ is an open cover for } I\}\geq \ent(G,\alpha_0)
$$
$$
=\lim_{m\to \infty}\frac{\log N(\star_{i=1}^{m}G^{-1},\alpha_0^{m+1})}{m}=\lim_{m\to \infty}\frac{\log N(\star_{i=1}^{i_1+i_2+i_3+\ldots +i_m+2m-1}G^{-1},\alpha_0^{i_1+i_2+i_3+\ldots +i_m+2m})}{i_1+i_2+i_3+\ldots +i_m+2m-1}
$$
$$
\geq  \lim_{m\to \infty}\frac{\log 2^m}{\psi^{L,R}_{b}\cdot m+2m}=\lim_{m\to \infty}\frac{m\log 2}{m\cdot (\psi^{L,R}_{b}+2)}=\frac{\log 2}{\psi^{L,R}_{b}+2}.
$$
Therefore, $\ent(G)>0$.
\end{proof}
The following example shows that Theorem \ref{MAIN} is not a characterization of the non-zero entropy for closed relations on $[0,1]$.
\begin{example}
Let $G=\{(0,1),(0,\frac{3}{4}),(\frac{3}{4},0),(1,0)\}$ be a closed relation on $[0,1]$.
It is easy to see that $\ent(G)\neq 0$ while there are no well-aligned  closed subsets of $G$ nor there are well-aligned  closed subsets of $G^{-1}$.
\end{example}
\section{Examples}\label{s5}
Finally, we present various examples of closed relations $G$ on $[0,1]$ such that
\begin{enumerate}
\item the entropy of $G$ is non-zero,
\item no periodic point or exactly one periodic point is generated by $G$, and
\item no Cantor set is finitely generated by $G$.
\end{enumerate}
In the following theorem, a closed relation on  $I$ which is  i-embedded into $I\times I$, is obtained.
\begin{theorem}\label{thm101}
Let $a\in (1,\infty)$  be any irrational number such that $a^k$ is irrational for any positive integer $k$ and let $b\in (0,1)$ be any rational number such that $\frac{1}{a}>b$.   Also, let
$$
H_{a,b}=\Big\{(x,y)\in I\times I \ | \ y=ax, x\in \Big[\frac{b}{a^2},\frac{1}{a}\Big] \Big\}\cup \Big\{(x,y)\in I\times I \ | \ y=bx, x\in \Big[\frac{1}{a^2},1\Big]\Big\}
$$
(the relation $H_{a,b}$ is pictured in Figure \ref{figure1}).
\begin{figure}[h!]
	\centering
		\includegraphics[width=20em]{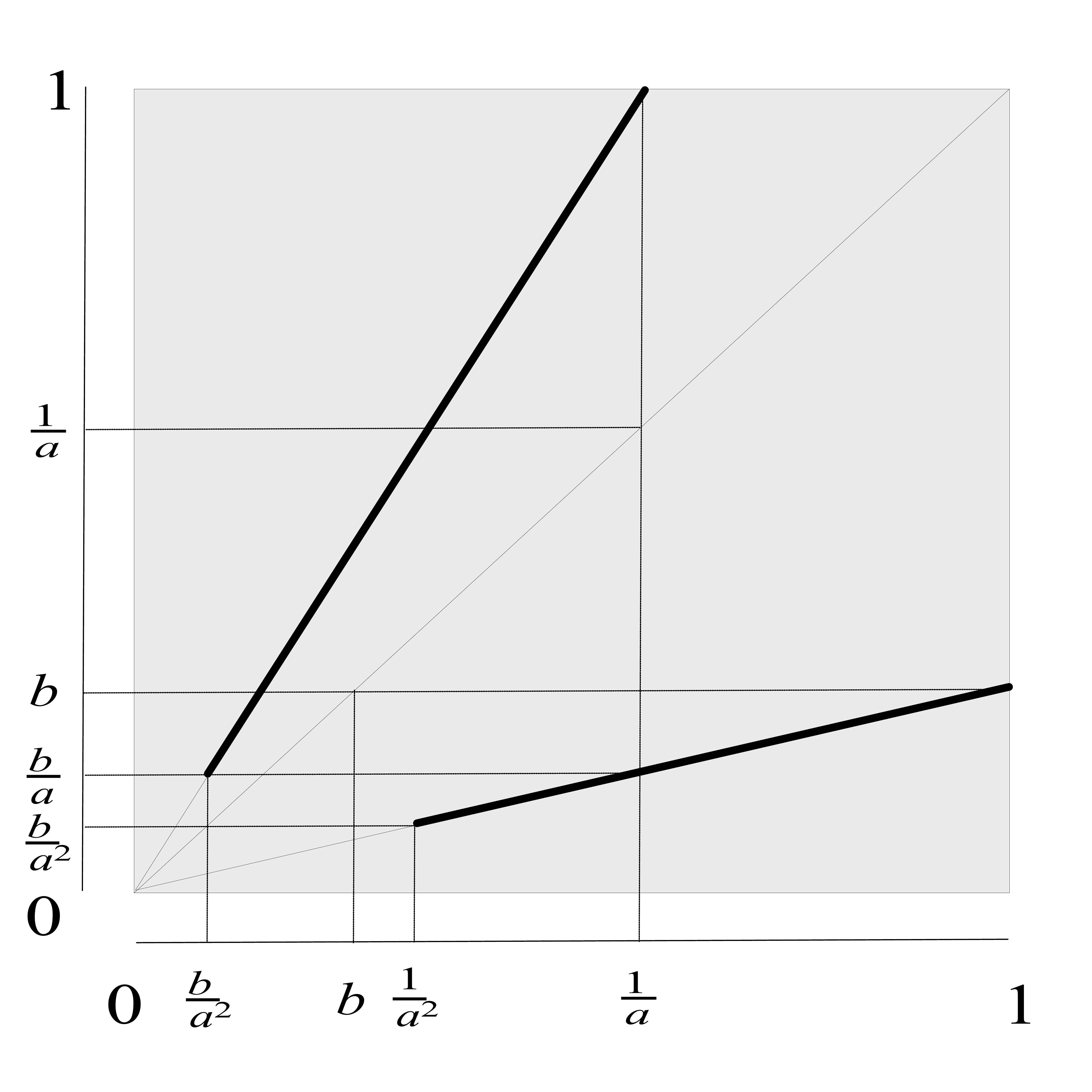}
	\caption{The relation $H_{a,b}$}
	\label{figure1}
\end{figure}  
Then $H_{a,b}$ is  i-embedded into $I\times I$.
\end{theorem}
\begin{proof} 
Taking
$$
L=\Big\{(x,y)\in I\times I \ | \ y=ax, x\in \Big[\frac{b}{a^2},\frac{1}{a}\Big] \Big\}
$$
and 
$$
R=\Big\{(x,y)\in I\times I \ | \ y=bx, x\in \Big[\frac{1}{a^2},1\Big]\Big\},
$$
it follows from Theorem \ref{MAIN}  that $\ent(H_{a,b})\neq 0$.

Next, we prove that $H_{a,b}$ does not generate a periodic point.  Note that for each positive integer $k$, $x_{k+1}=\frac{1}{a}x_k$ or  $x_{k+1}=\frac{1}{b}x_k$ holds. Suppose that there is a point $\mathbf x=(x_1,x_2,x_3,\ldots)\in \star_{i=1}^{\infty}H_{a,b}^{-1}$ and a positive integer $m$ such that $\sigma^m(\mathbf x)=\mathbf x$.  Then $x_1=x_{m+1}$ and, therefore, 
$$
x_{m+1}=\frac{1}{a^{k}}\cdot \frac{1}{b^{\ell}}\cdot x_1=x_1
$$
for some positive integers $k,\ell$. It follows that $\frac{1}{a^{k}}\cdot \frac{1}{b^{\ell}}=1$ and that 
$$
a^k=\frac{1}{b^{\ell}}
$$
meaning that $a^k$ is rational---a contradiction.  Therefore, $H_{a,b}$ does not generate a periodic point.
\end{proof}
Note that it also follows from Theorem \ref{gaga} that no Cantor set is finitely generated by $H_{a,b}$.

In the following theorem,  the graph $\Gamma(f)$ of an upper semi-continuous function  $f:I\multimap I$ such that $\Gamma(f)$ is  i-embedded into $I\times I$ is presented.

\begin{theorem}\label{thm11}
Let $a\in (1,\infty)$  be any irrational number such that $a^k$ is irrational for any positive integer $k$ and let $b\in (0,1)$ be any rational number such that $\frac{1}{a}>b$.   Also, let 
$$
H=H_{a,b}\cup\Big(\Big[0,\frac{b}{a^2}\Big]\times\Big\{\frac{b}{a}\Big\}\Big)\subseteq I\times I
$$
(the relation $H$ is pictured in Figure \ref{figure11}).  Then 
  \begin{enumerate}
\item $H$ is the graph of an upper semicontinuous function from $I$ to $I$, 
\item $H$ is  i-embedded into $I\times I$.
\end{enumerate}
\begin{figure}[h!]
	\centering
		\includegraphics[width=20em]{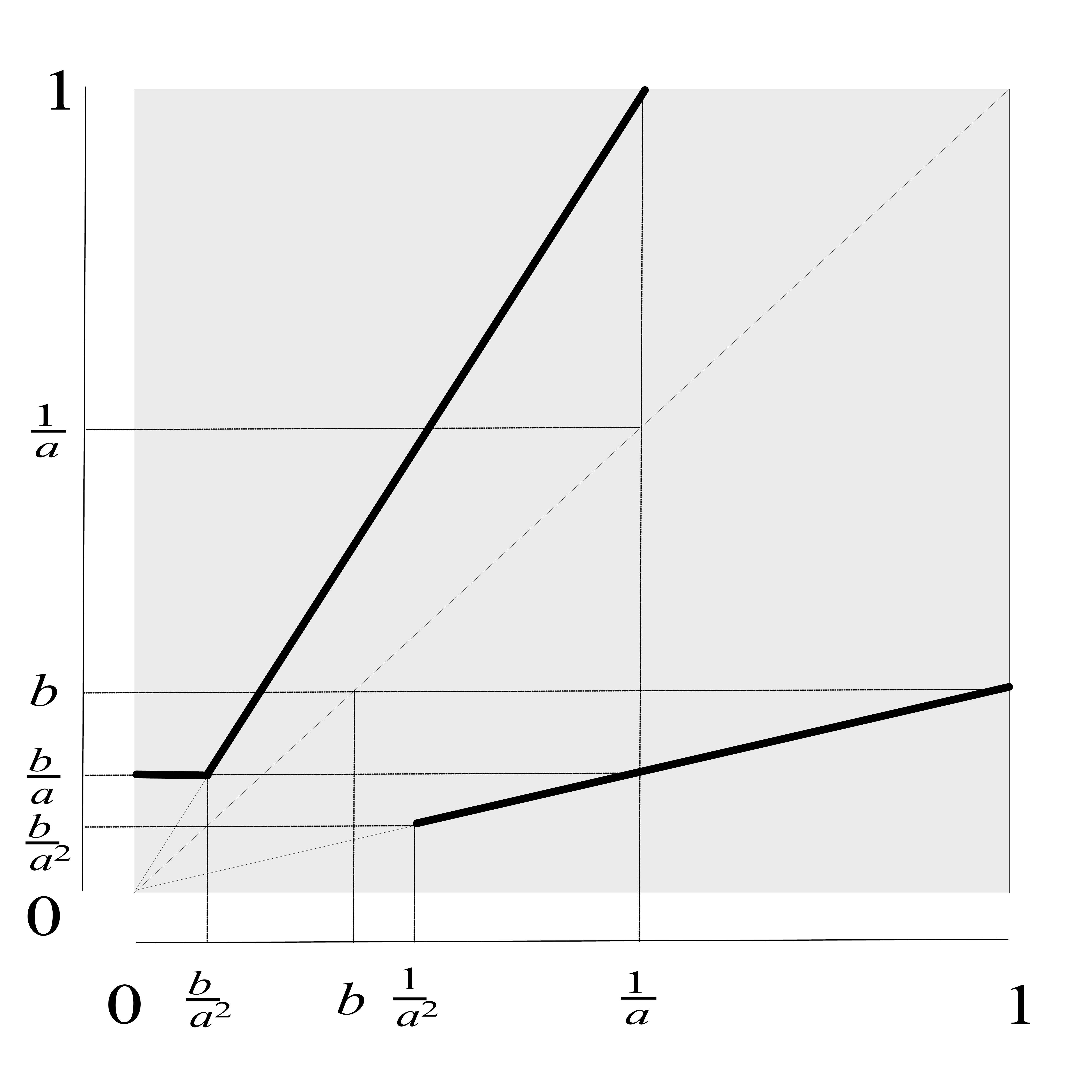}
	\caption{The relation $H$ from Theorem \ref{thm11}}
	\label{figure11}
\end{figure}  
\end{theorem}
\begin{proof}
First note that since $H$ is closed in $I\times I$ and since $\pi_1(H)=I$, it follows that $H$ is the graph of an upper semicontinuous function from $I$ to $I$. 

It follows from $H_{a,b}\subseteq H$ that $\ent(H)\neq 0$ by Theorem \ref{thm101}. 

Next, we prove that $H$ does not generate a periodic point.  First, note that for each point $\mathbf x=(x_1,x_2,x_3,\ldots)\in \star_{i=1}^{\infty}H^{-1}$, 
$$
\frac{b}{a^2}\leq x_1,x_2,x_3,\ldots \leq 1.
$$
Therefore, for each positive integer $k$, $x_{k+1}=\frac{1}{a}x_k$ or  $x_{k+1}=\frac{1}{b}x_k$ holds. Suppose that there is a point $\mathbf x=(x_1,x_2,x_3,\ldots)\in \star_{i=1}^{\infty}H^{-1}$ and a positive integer $m$ such that $\sigma^m(\mathbf x)=\mathbf x$.  Then $x_1=x_{m+1}$ and, therefore, 
$$
x_{m+1}=\frac{1}{a^{k}}\cdot \frac{1}{b^{\ell}}\cdot x_1=x_1
$$
for some positive integers $k,\ell$. It follows that $\frac{1}{a^{k}}\cdot \frac{1}{b^{\ell}}=1$ and that 
$$
a^k=\frac{1}{b^{\ell}}
$$
meaning that $a^k$ is rational---a contradiction.  Therefore, no periodic point is generated by $H$
\end{proof}
Note also that it follows from Theorem \ref{gaga},  that no Cantor set is finitely generated by $H$. 
\begin{observation}
Note that if $F:I\multimap I$ is an upper semi-continuous function such that  $\Gamma(F)$ is a subcontinuum of $I\times I$, then 
$$
\Gamma(F)\cap \{(t,t) \ | \ t\in [0,1]\}\neq \emptyset.
$$
Therefore, such a set $\Gamma(F)$ does generate a periodic point; i.e., if $(t_0,t_0)\in \Gamma(F)\cap \{(t,t) \ | \ t\in [0,1]\}$, then $(t_0,t_0,t_0,\ldots)$ is one of the periodic points,  generated by $\Gamma(F)$.
\end{observation}
The following theorem produces the best possible scenario for upper semi-continuous functions whose graphs are connected; i.e.,  the  graph $\Gamma(F)$ of a surjective  upper semi-continuous function  $F:I\multimap I$ such that $\Gamma(F)$ is a continuum and almost  i-embedded into $I\times I$, is presented.
\begin{theorem}\label{thm2}
Let $a\in (1,\infty)$  be any irrational number such that $a^k$ is irrational for any positive integer $k$ and let $b\in (0,1)$ be any rational number such that $\frac{1}{a}>b$.   Also, let
$$
H=\{(x,y)\in I\times I \ | \ y=ax \textup{ or } y=bx\}
$$
  (the relation $H$ is pictured in Figure \ref{figure2}). Then 
\begin{enumerate}
\item $H$ is the graph of a surjective  upper semicontinuous function from $I$ to $I$,
\item $H$ is a subcontinuum of $I\times I$, 
\item $H$ is almost  i-embedded into $I\times I$.
\end{enumerate}
\begin{figure}[h!]
	\centering
		\includegraphics[width=20em]{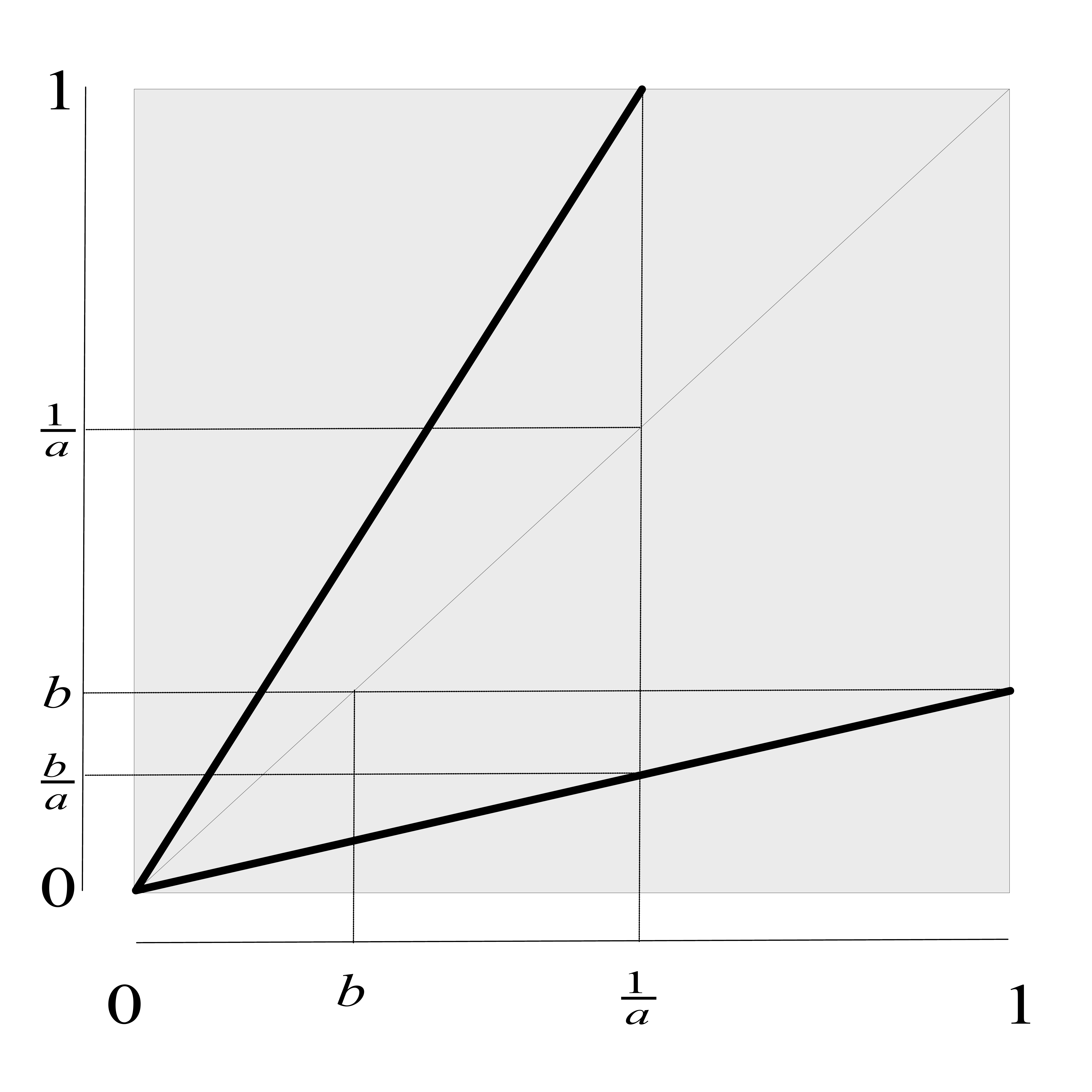}
	\caption{The relation $H$ from Theorem \ref{thm2}}
	\label{figure2}
\end{figure}
\end{theorem}
\begin{proof}

First, note that since $H$ is closed in $I\times I$, and since $p_1(H)=I$ and $p_2(H)=I$, it follows that $H$ is the graph of an upper semicontinuous function from $I$ to $I$ whose graph is surjective.  It is obvious that $H$ is a subcontinuum of $I\times I$. 

It follows from $H_{a,b}^{-1}\subseteq H$ that $\ent(H)\neq 0$ by Theorem \ref{thm101}.

Next, we prove that $\mathbf o=(0,0,0,\ldots )$ is the only periodic point, generated by $H$.  First, note that for each point $\mathbf x=(x_1,x_2,x_3,\ldots)\in \star_{i=1}^{\infty}H^{-1}$,  $x_{k+1}=\frac{1}{a}x_k$ or  $x_{k+1}=\frac{1}{b}x_k$ holds for each positive integer $k$. 
Suppose that there is a point $\mathbf x=(x_1,x_2,x_3,\ldots)\in \star_{i=1}^{\infty}H^{-1}\setminus\{\mathbf o\}$ and a positive integer $m$ such that $\sigma^m(\mathbf x)=\mathbf x$ (note that $x_k\neq 0$ for each positive integer $k$). Then $x_1=x_{m+1}$ and, therefore, 
$$
x_{m+1}=\frac{1}{a^{k}}\cdot \frac{1}{b^{\ell}}\cdot x_1=x_1
$$
for some positive integers $k,\ell$. It follows that $\frac{1}{a^{k}}\cdot \frac{1}{b^{\ell}}=1$ and that 
$$
a^k=\frac{1}{b^{\ell}}
$$
meaning that $a^k$ is rational---a contradiction.  Therefore, $H$ does not generate any other periodic points.
\end{proof}
Note that by Theorem \ref{gaga},  no Cantor set is finitely generated by $H$. 
\begin{lemma}\label{joj3}
Let $a\in (1,\sqrt{2})$. Then $\frac{a-1}{a+1}<\frac{1}{a}$.
\end{lemma}
\begin{proof}
The lemma follows from the observation that for each $a\in \mathbb R\setminus\{0,-1\}$,
$$
\frac{a-1}{a+1}<\frac{1}{a} \Longleftrightarrow a^2-2a-1<0 \Longleftrightarrow a\in (1-\sqrt{2},0)\cup(0,1+\sqrt{2}).  
$$
\end{proof}
\begin{lemma}\label{joj4}
Let $a\in (1,\sqrt{2})$. Then $\frac{a}{a+1}<\frac{1+a}{2a}$.
\end{lemma}
\begin{proof}
The lemma follows from the observation that for each $a\in \mathbb R\setminus\{0,-1\}$,
$$
\frac{a}{a+1}<\frac{1+a}{2a} \Longleftrightarrow a^2-2a-1<0 \Longleftrightarrow a\in (1-\sqrt{2},0)\cup(0,1+\sqrt{2}).
$$
\end{proof}
\begin{observation}
Let $a\in (1,\sqrt{2})$ and let $b\in (\frac{a}{a+1},\frac{1+a}{2a})$.  Then $a>1$ and $b<1$, and, therefore, 
$$
\frac{2a}{1+a}\in (1,\infty)
$$
and 
$$
\frac{b+1}{2}\in(0,1).
$$
\end{observation}
\begin{lemma}\label{joj5}
Let $a\in (1,\sqrt{2})$ and $b\in (\frac{a}{a+1},\frac{1+a}{2a})$, let $A$ be the closed relation on $I$, defined by
$$
A=\Big\{(x,y)\in [0,\frac{a+1}{2a}]\times [0,1] \ | \ y=\frac{2a}{1+a}x \Big\}
$$
$$
\cup \Big\{(x,y)\in [0,1]\times [0,\frac{b+1}{2}] \ | \  y=\frac{b+1}{2}x\Big\},
$$
and let $B$ be the closed relation on $[-1,1]$,  defined by
$$
B=\Big\{(x,y)\in [-1,\frac{1}{a}]\times [-1,1] \ | \ y=\frac{2a}{a+1}x +\frac{a-1}{a+1}\Big\}
$$
$$
\cup \Big\{(x,y)\in [-1,1]\times [-1,b] \ | \  y=\frac{b+1}{2}x+\frac{b-1}{2}\Big\}.
$$
Then $A$ and $B$ are topological conjugates.
\end{lemma}
\begin{proof}
Let
$$
\varphi:[-1,1]\rightarrow [0,1]
$$
be the homeomorphism, defined by
$$
\varphi(t)=\frac{1}{2}t+\frac{1}{2}
$$
for each $t\in [-1,1]$. We show that for each $(s,t)\in [-1,1]\times [-1,1]$, 
 $$
 (s,t)\in B  \Longleftrightarrow (\varphi(s), \varphi(t))\in A.
 $$
Let $(s,t)\in B$ be any point.  We treat the following possible cases.
\begin{enumerate}
\item $(s,t)\in \{(x,y)\in [-1,\frac{1}{a}]\times [-1,1] \ | \ y=\frac{2a}{a+1}x +\frac{a-1}{a+1}\}$. Then $t=\frac{2a}{a+1}s +\frac{a-1}{a+1}$ and, therefore,
$$
(\varphi(s),\varphi(t))=\Big(\frac{1}{2}s+\frac{1}{2},\frac{a}{a+1}s+\frac{a}{a+1}\Big)
$$
while 
$$
\frac{2a}{a+1}\varphi(s)=\frac{2a}{a+1}\Big(\frac{1}{2}s+\frac{1}{2}\Big)=\frac{a}{a+1}s+\frac{a}{a+1}.
$$
Since $\varphi(s)\in [0,\frac{a+1}{2a}]$, it follows that 
$$
(\varphi(s), \varphi(t))\in \Big\{(x,y)\in \Big[0,\frac{a+1}{2a}\Big]\times [0,1] \ | \ y=\frac{2a}{1+a}x \Big\}
$$
 and $(\varphi(s), \varphi(t))\in A$ follows. 
\item $(s,t)\in  \{(x,y)\in [-1,1]\times [-1,b] \ | \  y=\frac{b+1}{2}x+\frac{b-1}{2}\}$.  Then $t=\frac{b+1}{2}s+\frac{b-1}{2}$ and, therefore,
$$
(\varphi(s),\varphi(t))=\Big(\frac{1}{2}s+\frac{1}{2},\frac{b+1}{4}s+\frac{b+1}{4}\Big)
$$
while 
$$
\frac{b+1}{2}\varphi(s)=\frac{b+1}{2}\Big(\frac{1}{2}s+\frac{1}{2}\Big)=\frac{b+1}{4}s+\frac{b+1}{4}.
$$
Since $\varphi(t)\in [0,\frac{b+1}{2}]$, it follows that 
$$
(\varphi(s), \varphi(t))\in \Big\{(x,y)\in [0,1]\times \Big[0,\frac{b+1}{2}\Big] \ | \  y=\frac{b+1}{2}x\Big\}
$$
 and $(\varphi(s), \varphi(t))\in A$ follows. 
\end{enumerate}
Next, let $(\varphi(s), \varphi(t))\in A$.  To show that $(s,t)\in B$, we treat the following possible cases.
\begin{enumerate}
\item $(\varphi(s), \varphi(t))\in \{(x,y)\in [0,\frac{a+1}{2a}]\times [0,1] \ | \ y=\frac{2a}{1+a}x \}$. Then $\varphi(t)=\frac{2a}{a+1}\varphi(s)$.  It follows that 
$$
\frac{1}{2}t+\frac{1}{2}=\frac{2a}{a+1}\Big(\frac{1}{2}s+\frac{1}{2}\Big)=\frac{a}{a+1}s+\frac{a}{a+1}
$$
and, therefore,
$$
t=\frac{2a}{a+1}s +\frac{a-1}{a+1}.
$$
Since $s\in [-1,\frac{1}{a}]$, it follows that 
$$
(s, t)\in \Big\{(x,y)\in \Big[-1,\frac{1}{a}\Big]\times [-1,1] \ | \ y=\frac{2a}{a+1}x +\frac{a-1}{a+1}\Big\}
$$
 and $(s,t)\in B$ follows. 
\item $(\varphi(s), \varphi(t))\in \{(x,y)\in [0,1]\times [0,\frac{b+1}{2}] \ | \  y=\frac{b+1}{2}x\}$. Then $\varphi(t)=\frac{b+1}{2}\varphi(s)$.  It follows that 
$$
\frac{1}{2}t+\frac{1}{2}=\frac{b+1}{2}\Big(\frac{1}{2}s+\frac{1}{2}\Big)=\frac{b+1}{4}s+\frac{b+1}{4}
$$
and, therefore,
$$
t=\frac{b+1}{2}s +\frac{b-1}{2}.
$$
Since $t\in [-1,b]$, it follows that 
$$
(s, t)\in\Big\{(x,y)\in [-1,1]\times [-1,b] \ | \  y=\frac{b+1}{2}x+\frac{b-1}{2}\Big\}
$$
 and $(s,t)\in B$ follows. 
\end{enumerate}
\end{proof}
\begin{lemma}\label{ffff}
Let $a\in (1,\sqrt{2})$ and $b\in (\frac{a}{a+1},\frac{1+a}{2a})$, and let $H$ be the closed relation on $I$ that is defined by
$$
H=\Big\{(x,y)\in [0,\frac{1}{a}]\times [0,1] \ | \ y=\frac{2a}{a+1}x +\frac{a-1}{a+1}\Big\}
$$
$$
\cup \Big\{(x,y)\in [0,1]\times [0,b] \ | \  y=\frac{b+1}{2}x+\frac{b-1}{2}\Big\}
$$
(the relation $H$ is pictured in Figure \ref{tatatale}). 
\begin{figure}[h!]
	\centering
		\includegraphics[width=20em]{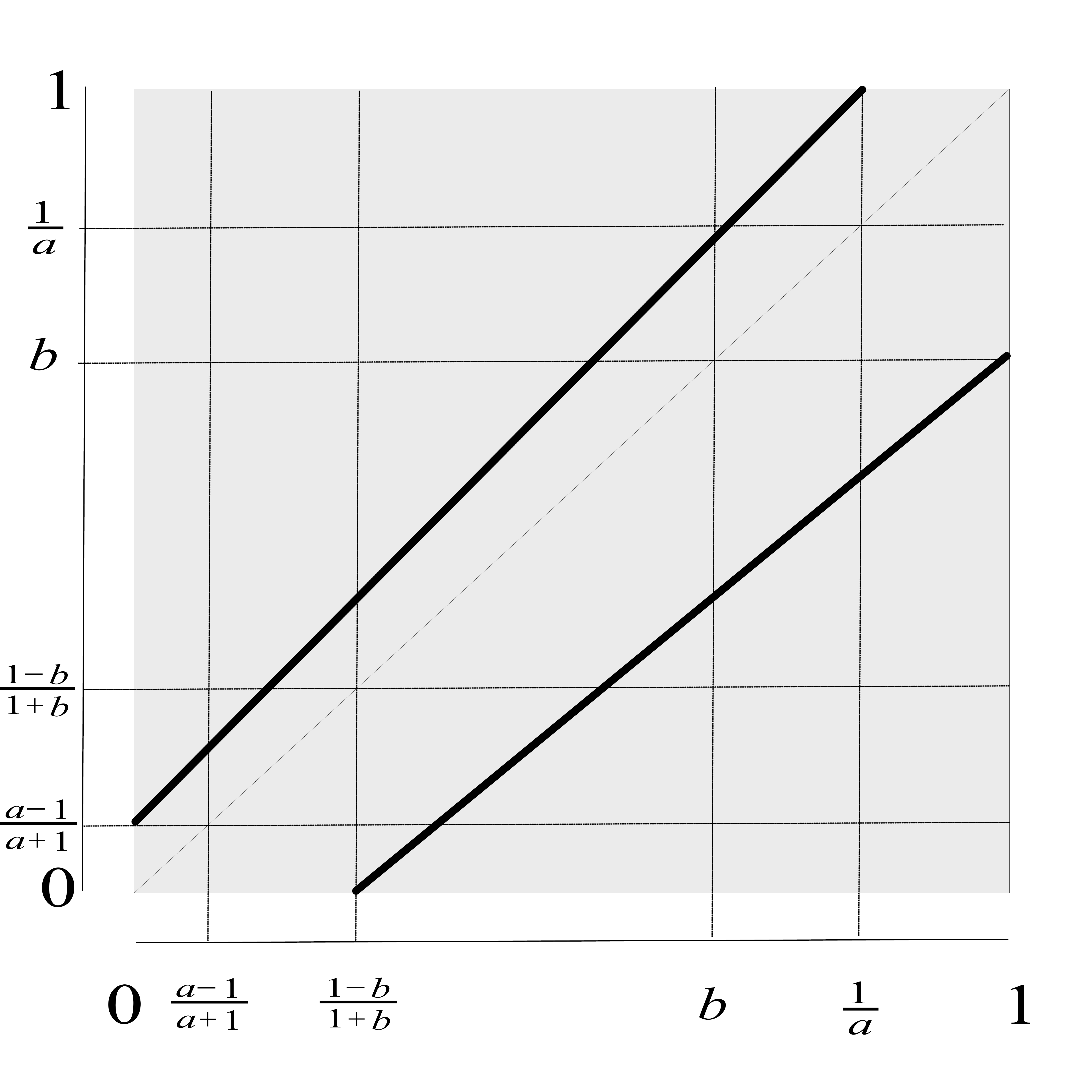}
	\caption{The relation $H$ from Lemma \ref{ffff}}
	\label{tatatale}
\end{figure}  
Then  there is a surjective upper semicontinuous function $F:I\multimap I$ such that $\Gamma(F)=H$. 
\end{lemma}
\begin{proof}
Obviously, $H$ is closed in $I\times I$.  The set $\{(x,y)\in [0,\frac{1}{a}]\times [0,1] \ | \ y=\frac{2a}{a+1}x +\frac{a-1}{a+1}\}$ is the straight line segment in $I\times I$ with endpoints $(0,\frac{a-1}{a+1})$ and $(\frac{1}{a},1)$ and the set $\{(x,y)\in [0,1]\times [0,b] \ | \  y=\frac{b+1}{2}x+\frac{b-1}{2}\}$ is the straight line segment in $I\times I$ with endpoints $(\frac{1-b}{1+b},0)$ and $(1,b)$.  

Note that 
$$
b>\frac{a-1}{a+1} \textup{ and } \frac{1-b}{1+b}<\frac{1}{a}
$$
since $b\in (\frac{a}{a+1},\frac{1+a}{2a})$. Therefore, $p_1(H)=I$ and $p_2(H)=I$.  This means that there is a surjective upper semicontinuous function $F:I\multimap I$ such that $\Gamma(F)=H$. 
\end{proof}
\begin{lemma}\label{jojmene}
Let $a\in (1,\sqrt{2})$ such that $\left(\frac{2a}{1+a}\right)^k$
is not rational for any positive integer $k$, and let $b\in (\frac{a}{a+1},\frac{1+a}{2a})$ be a rational number.  Also, let $H$ be the closed relation on $I$ defined by
$$
H=\Big\{(x,y)\in [0,\frac{1}{a}]\times [0,1] \ | \ y=\frac{2a}{a+1}x +\frac{a-1}{a+1}\Big\}
$$
$$
\cup \Big\{(x,y)\in [0,1]\times [0,b] \ | \  y=\frac{b+1}{2}x+\frac{b-1}{2}\Big\}
$$
(the relation $H$ is pictured in Figure \ref{tatatale}).  Then  no periodic point is generated by $H$.
\end{lemma}
\begin{proof}
Let $A$ be the closed relation on $[0,1]$ defined by
$$
A=\Big\{(x,y)\in [0,\frac{a+1}{2a}]\times [0,1] \ | \ y=\frac{2a}{1+a}x \Big\}
$$
$$
\cup \Big\{(x,y)\in [0,1]\times [0,\frac{b+1}{2}] \ | \  y=\frac{b+1}{2}x\Big\},
$$
and let $B$ be the closed relation on $[-1,1]$ defined by
$$
B=\Big\{(x,y)\in [-1,\frac{1}{a}]\times [-1,1] \ | \ y=\frac{2a}{a+1}x +\frac{a-1}{a+1}\Big\}
$$
$$
\cup \Big\{(x,y)\in [-1,1]\times [-1,b] \ | \  y=\frac{b+1}{2}x+\frac{b-1}{2}\Big\}.
$$
The relations  $A$ and $B$ are topological conjugates by Lemma \ref{joj5}.
Let $\varphi:[-1,1]\rightarrow [0,1]$ be the homeomorphism, defined by $\varphi(t)=\frac{1}{2}t+\frac{1}{2}$ for each $t\in [-1,1]$.  It follows from Theorem \ref{thm2}  that $\mathbf o=(0,0,0,\ldots )$ is the only periodic point in $\star_{i=1}^{\infty}A^{-1}$, generated by $ A$. Therefore, by Lemma \ref{joj0},  $(-1,-1,-1,\ldots)$ is the only periodic point in $\star_{i=1}^{\infty}B^{-1}$, generated by $ B$.  Since $(-1,-1,-1,\ldots)\not \in \star_{i=1}^{\infty}H^{-1}$ and since $H\subseteq B$,  no periodic point is generated by $H$. 
\end{proof}

\begin{observation}\label{loll}
Let $a\in (1,\sqrt{2})$ and $b\in (\frac{a}{a+1},\frac{1+a}{2a})$.
Note that
\begin{enumerate}
\item $\frac{1-b}{1+b}>\frac{a-1}{a+1}$ if and only if $\frac{1}{a}>b$.
\item $\frac{1}{a}<\frac{1+a}{2a}$ if and only if $a(a-1)>0$, and this is true for any $a\in (1,\sqrt{2})$.
\item $\frac{a}{a+1}<\frac{1}{a}$ if and only if $a\in (\frac{1-\sqrt{5}}{2},\frac{1+\sqrt{5}}{2})\setminus\{-1,0\}$,  and since  $\sqrt{2}<\frac{1+\sqrt 5}{2}$,  
$$
\frac{a}{a+1}<\frac{1}{a}
$$
 is true for any $a\in (1,\sqrt{2})$.
\end{enumerate}
So, for any $a\in (1,\sqrt{2})$, 
$$
\frac{1}{a}\in (\frac{a}{a+1},\frac{1+a}{2a}).
$$
To achieve that $\frac{1}{a}>b$, we will always take $b\in (\frac{a}{a+1},\frac{1}{a})$. This way
$$
\frac{1-b}{1+b}>\frac{a-1}{a+1}
$$
is achieved for any  $a\in (1,\sqrt{2})$ and $b\in (\frac{a}{a+1},\frac{1}{a})$.
\end{observation}

Next, we present a theorem, where the  graph $\Gamma(F)$ of a surjective  upper semi-continuous function  $F:I\multimap I$ such that $\Gamma(F)$ is  i-embedded into $I\times I$.

\begin{figure}[h!]
	\centering
		\includegraphics[width=20em]{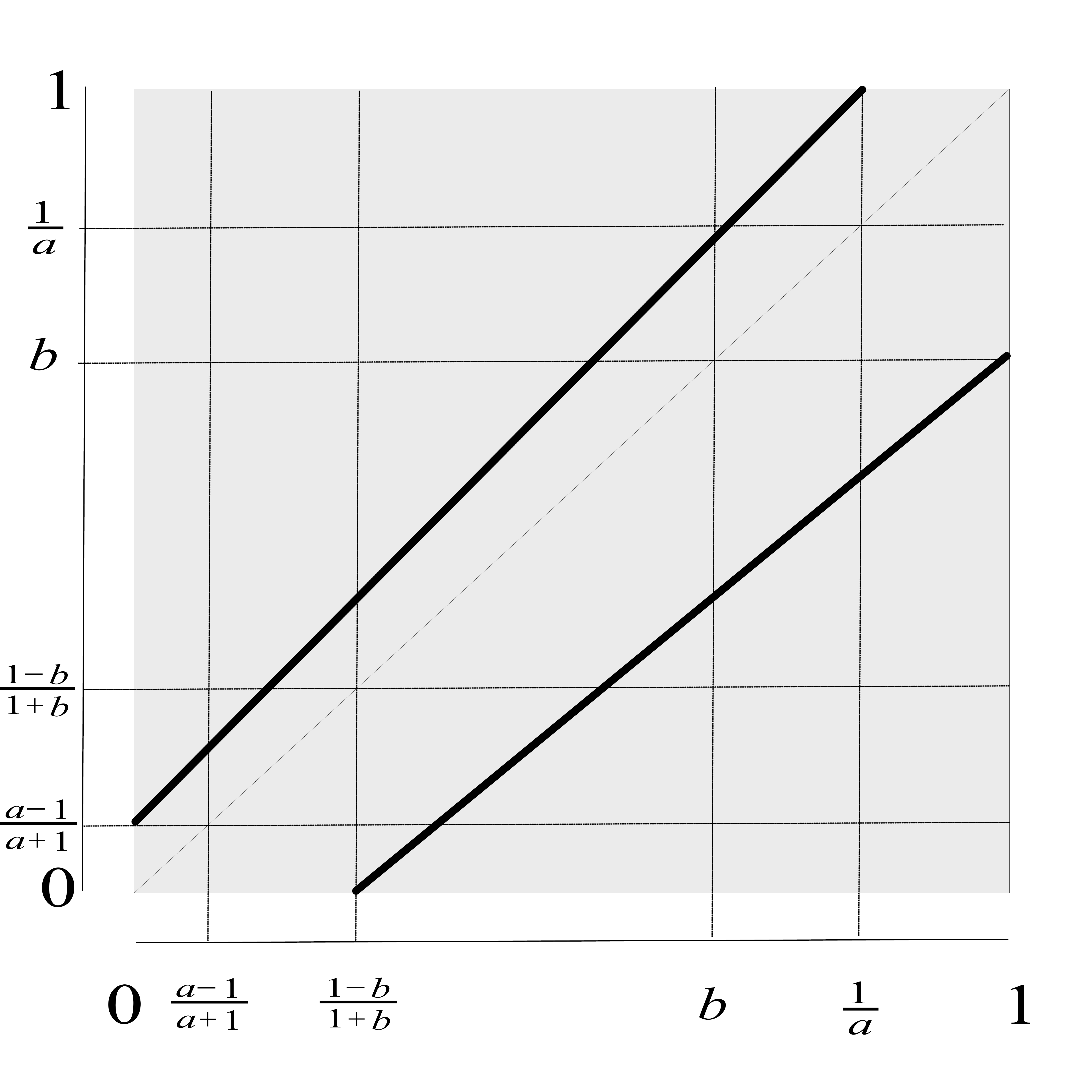}
	\caption{The relation $H$ from Theorem \ref{taletoti}}
	\label{tatatale11}
\end{figure} 

\begin{theorem}\label{taletoti}
Let $a\in (1,\sqrt{2})$ such that $\left(\frac{2a}{1+a}\right)^k$
is not rational for any positive integer $k$, and let $b\in (\frac{a}{a+1},\frac{1}{a})$ be a rational number,  and let $H$ be the closed relation on $I$, defined by
$$
H=\Big\{(x,y)\in [0,\frac{1}{a}]\times [0,1] \ | \ y=\frac{2a}{a+1}x +\frac{a-1}{a+1}\Big\}
$$
$$
\cup \Big\{(x,y)\in [0,1]\times [0,b] \ | \  y=\frac{b+1}{2}x+\frac{b-1}{2}\Big\}
$$
(the relation $H$ is pictured in Figure \ref{tatatale11}).  
Then 
\begin{enumerate}
\item $H$ is the graph of a surjective  upper semicontinuous function from $I$ to $I$,
\item $H$ is  i-embedded into $I\times I$.
\end{enumerate} 
\end{theorem}
\begin{proof}
By Lemma \ref{ffff}, $H$ is the graph of a surjective  upper semicontinuous function from $I$ to $I$.  By Lemma \ref{jojmene}, no periodic point is generated by $H$. 

Finally,  by letting 
$$
L=\Big\{(x,y)\in [0,\frac{1}{a}]\times [0,1] \ | \ y=\frac{2a}{a+1}x +\frac{a-1}{a+1}\Big\}
$$
and
$$
R=\Big\{(x,y)\in [0,1]\times [0,b] \ | \  y=\frac{b+1}{2}x+\frac{b-1}{2}\Big\},
$$
it follows from Theorem \ref{MAIN} that the entropy of $H$ is non-zero.  Therefore,  $H$ is  i-embedded into $I\times I$.
\end{proof}



\noindent I. Bani\v c\\
              (1) Faculty of Natural Sciences and Mathematics, University of Maribor, Koro\v{s}ka 160, SI-2000 Maribor,
   Slovenia; \\(2) Institute of Mathematics, Physics and Mechanics, Jadranska 19, SI-1000 Ljubljana, 
   Slovenia; \\(3) Andrej Maru\v si\v c Institute, University of Primorska, Muzejski trg 2, SI-6000 Koper,
   Slovenia\\
             {iztok.banic@um.si}           
     
				\-
				
		\noindent G.  Erceg\\
             Faculty of Science, University of Split, Rudera Bo\v skovi\' ca 33, Split,  Croatia\\
{{goran.erceg@pmfst.hr}       }    

                 	\-
					
  \noindent J.  Kennedy\\
             Department of Mathematics, Lamar University, 200 Lucas Building, P.O. Box 10047, Beaumont, TX 77710 USA\\
{{kennedy9905@gmail.com}       }    


\begin{thebibliography}{9}


\bibitem[AKM] {AKM} R. L. Adler, G. Konheim, M. H. McAndrew, Topological entropy \textit{Transactions of the American Mathematical Society} \textbf{114} (1965) 309-319.

\bibitem[A] {A} E. Akin, \textbf{General Topology of Dynamical Systems}, Volume 1, Graduate Studies in Mathematics Series, American Mathematical Society, Providence RI, 1993.

\bibitem[B]{B} Iztok Bani\v{c}, Inverse limits as limits with respect to the
Hausdorff metric, \textit{Bulletin of the Australian Mathematical Society} 
\textbf{75} (2007) 17-22.


\bibitem[BCMM1]{BCMM1}Iztok Bani\v{c}, Matev\v{z} \v{C}repnjak, Matej Merhar, Uro
\v{s} Milutinovi\v{c}, Limits of inverse limits, \textit{Topology Appl.}, \textbf{157} (2010) 439-450.

\bibitem[BCMM2] {BCMM2}  I. Bani\v{c}, M. \v{C}repnjak. M. Merhar, U. Milutinovi\v{c}, Towards the complete classification of tent maps  inverse limits, \textit{Topology and its Applications} \textbf{160} (2013) 63-73.

\bibitem[BCMM3] {BCMM3} I. Bani\v{c}, M. \v{C}repnjak. M. Merhar, U. Milutinovi\v{c}, T. Sovi\v{c},  Wa\.{z}ewski's universal dendrite as an inverse limit with one set-valued bonding function, \textit{Glas. Mat.} \textbf{48} (2013) 137-165.

\bibitem[BK] {BK} I. Bani\v{c} and J. Kennedy, Inverse limits with bonding functions whose graphs are arcs, \textit{Topology and its Applications} \textbf{151} (2015) 9 - 21.

\bibitem[BBCDS]{BBCDS} J. Banks, J. Brooks, G. Cairns, G. Davis, P. Stacey, On Devaney's Definition of Chaos, \textit{The American Mathematical Monthly} \textbf{99}, No. 4 (1992), pp. 332-334

\bibitem[Bo1] {Bo} R. E. Bowen, Entropy for group endomorphisms and homogeneous spaces, \textit{Transactions of the American Mathematical Society} \textbf{153} (1971) 401-414.

\bibitem[Bo2] {Bo2} R. E. Bowen, Topological entropy and axiom A, \textit{Global Analysis (Proc. Sympos. Pure Math., Vol. XIV, Berkeley, Calif., 1968), Amer. Math. Soc.}, Providence, R.I., 1970, pp. 23--41.

\bibitem[CC]{CC} J. H. Case and R. E. Chamberlin, Characterizations of tree-like continua, \textit{Pacific J. Mathematics} \textbf{10} (1960) 73-84.

\bibitem[CR]{CR} W. J. Charatonik and R. P. Roe, Inverse limits of continua having trivial shape, \textit{Houston J.l of Mathematics} \textbf{38}, no. 4, (2012) 1307-1312.

\bibitem[D] {D} R. Devaney, \textbf{ An Introduction to Chaotic Dynamical Systems, 2nd. ed.}, Westview Press, Cambridge MA, 2003.

\bibitem[Di] {Di} E. I. Dinaburg, The relation between topological entropy and metric entropy, \textit{Soviet Math.} \textbf{11} (1970) 13-16.

\bibitem[EK] {EK} G. Erceg and J. Kennedy, Topological entropy on closed sets in $[0,1]^2$, \textit{Topology Appl.} \textbf{246} (2018) 106-136. 


\bibitem[GK1]{GK1} Sina Greenwood and Judy Kennedy, Connected generalized inverse limits, \textit{Topology Appl.}, \textbf{159} (2012) 57-68.

\bibitem[GK2]{GK2} Sina Greenwood and Judy Kennedy, Connectedness and Ingram-Mahavier products, \textit{Topology Appl.}, \textbf{166} (2014) 1-9.

\bibitem[Il] {Il} A. Illanes, A circle is not the generalized inverse limit of a subset of $[0,1]^{2}$, \textit{Proceedings of the AMS} \textbf{139} (2011) 2987 -2993.

\bibitem[I1]{I1}W. T. Ingram, Two-pass maps and indecomposability of inverse
limits of graphs, \textit{Topology Proceedings} \textbf{29} (2005) 1-9.

\bibitem[I2]{I2}W. T. Ingram, Inverse limits of upper semicontinuous functions
that are unions of mappings, \textit{Topology Proceedings} \textbf{34}
(2009) 17-26.

\bibitem[I4]{I4}W. T. Ingram, Inverse limits of upper semicontinuous set valued functions,  \textit{Houston Journal of Mathematics} \textbf{32}, No 1,
(2006) 17-26.

\bibitem[I3]{I3} W. T. Ingram, Inverse limits with upper semicontinuous bonding
functions: Problems and some partial solutions, \textit{Topology Proc.}, \textbf{36} (2010) 353-373.

\bibitem[I5]{I5}W. T. Ingram, Two-pass maps and indecomposability of inverse
limits of graphs, \textit{Topology Proceedings} \textbf{29} (2005) 1-9.


\bibitem[I6] {I6} W. T. Ingram, \textbf{An Introduction to Inverse Limits with Set-Valued Functions}, Springer, New York, 2012.

  

\bibitem[I7]{I7} W. T.. Ingram, Inverse limits of upper semi-continuous functions that are unions of mappings, \textit{Topology Proceedings} \textbf{34} (2009) 17-26.

\bibitem[IM1]{IM1} W. T. Ingram and W. Mahavier, \textbf{Inverse Limits: From Continua to Chaos}, Springer, New York, NY, 2012.

 
\bibitem[IM2]{IM2}W. T. Ingram and William S. Mahavier, Inverse limits of upper
semicontinuous set valued functions, \textit{Houston Journal of Mathematics} \textbf{32} (2006) 119-130.

\bibitem[KT]{KT} J. Kelly and T. Tennant, Topological entropy of set-valued functions, \textit{Houston Journal of Mathematics} \textbf{43} (2015) 263-282.

\bibitem[KN]{KN} J. Kennedy and V. Nall, Dynamical properties of inverse limits with set valued functions, \textit{Ergodic Theory and Dynamical Systems}, publishied online Sept. 22, 2016.



\bibitem[L]{L} M. Lockyer, Dissertation: Topics in Generalised Inverse Limits, The University of Auckland, 2014.

 

\bibitem[M]{M}William S. Mahavier, Inverse limits with subsets of $[0,1]\times
\lbrack 0,1]$, \textit{Topology and its Applications} \textbf{141} (2004)
225-231.




\bibitem[Mi]{Mi} J. Mioduszewski, On a quasi-ordering in the class of continuous mappings of a closed interval into itself, \textit{Colloquium Mathematicum} \textbf{9} (1962) 233 - 240.

\bibitem[N1]{N1} Van Nall, Connected Inverse limits with set-valued functions,  \textit{Topology Proceedings} \textbf{40} (2012), 167-177.

\bibitem[N2]{N2} V. Nall, Inverse limits with set valued functions, \textit{Houston J. of Mathematics} \textbf{37}. no. 4 (2011) 1323-1332.

\bibitem[N3]{N3} V. Nall, Finite graphs that are inverse limits with a set valued function on $[0,1]$, \textit{Topology and its Applications} \textbf{158} (2011) 1226-1233.
 
\bibitem[N4]{N4} V. Nall,  The only finite graph that is an inverse limit with a set valued function on $[0,1]$ is an arc, \textit{Topology and its Applications} \textbf{159} (2012) 733-736.

\bibitem[RT]{RT} B. Raines and T. Tennant, The specification property on a set-valued map and its inverse limit, arxiv:1509.08415v1


\bibitem[R] {R} Clark Robinson, \textbf{Dynamical Systems: Stability, Symbolic Dynamics, and Chaos, 2nd ed.}, CRC Press, Boca Raton, 1998.


\bibitem[V]{V} S. Varagona, Inverse limits with upper semi-continuous bonding functions and indecomposability, \textit{Houston J. of Mathematics} \textbf{37} (2011) 1017-1034.

\bibitem[W] {Walters book} Peter Walters, \textbf{An Introduction to Ergodic Theory}, Springer-Verlag, New York, 1982.



\end{thebibliography}
\end{document}